\documentclass{amsart}
\usepackage[]{hyperref} 
\usepackage{float}
\usepackage{etoolbox}
\usepackage{color}
\patchcmd{\subsection}{-.5em}{.5em}{}{}
\newtheorem{theorem}{Theorem}[section]

\theoremstyle{definition}

\theoremstyle{remark}

\theoremstyle{prop}

\theoremstyle{coro}
\newtheorem{corollary}[theorem]{Corollary}
\numberwithin{equation}{section}

\makeatletter
\newcommand*{\rom}[1]{\expandafter\@slowromancap\romannumeral #1@}
\makeatother

\allowdisplaybreaks[1]
\begin{document}
\title[A Generic Classification of Exceptional Orthogonal X$_{1}$-Polynomials]{A Generic Classification of Exceptional Orthogonal X$_{1}$-Polynomials Based on Pearson Distributions Family}

\author[M. Masjed-Jamei]{Mohammad Masjed-Jamei$ ^{\dagger} $}

\author[Z. Moalemi]{Zahra Moalemi$ ^{\ddagger} $}

\thanks{ $ ^{\dagger},^{\ddagger} $ Department of Mathematics, K.N.Toosi University of Technology, P.O.Box 16315-1618, Tehran, Iran.\\
{\em $ ^{\dagger} $ E-mail address}, M. Masjed-Jamei: mmjamei@kntu.ac.ir ; mmjamei@yahoo.com\\
URL: https://wp.kntu.ac.ir/mmjamei\\
{\em $ ^{\ddagger}$ E-mail address}, Z. Moalemi: zmoalemi@mail.kntu.ac.ir}


\date{}

\begin{abstract}
The so-called exceptional orthogonal X$ _{1} $-polynomials arise as eigenfunctions of a Sturm-Liouville problem. In this paper, a generic classification of these polynomials is presented based on Pearson distributions family. Then, six special differential equations of the aforesaid classification are introduced and their polynomial solutions are studied in detail.
\end{abstract}


\maketitle

\section{Introduction}

Classical orthogonal polynomials are known to play a fundamental role in the construction of bound-state solutions to exactly solvable potentials in quantum mechanics \cite{MR2439200}. They are eigenfunctions of some Sturm-Liouville problems and form complete sets with respect to some positive-definite measures \cite{MR2191786,MR2656096,MR1149380}.

Consider the second order differential equation
\begin{equation}\label{EQ1.ODE}
\frac{d}{{dx}}\left( {k(x)\frac{{dy}}{{dx}}} \right) -\big(\lambda \rho(x)+ q(x)\big)y=0,
\end{equation}
on an open interval, say $ (a,b) $, with the boundary conditions
\begin{equation}\label{EQ2.boundary}
\begin{array}{l}
{\alpha _1}y(a)\, + {\beta _1}y'(a) = 0,\\[3mm]
{\alpha _2}y(b)\, + {\beta _2}y'(b) = 0,
\end{array}
\end{equation}
in which $ \alpha_{1},\alpha_{2} $ and $ \beta_{1},\beta_{2} $ are given constants and the functions $ k(x)>0\,,\,k'(x)\,,\,q(x) $ and $ \rho(x)>0 $ in \eqref{EQ1.ODE} are assumed  to be continuous for $ x\in [a,b] $. The boundary value problem \eqref{EQ1.ODE}-\eqref{EQ2.boundary} is called a regular Sturm-Lioville problem and if one of the points $ a $ and $ b $ is singular (i.e. $ k(a)=0 $ or $ k(b)=0 $), it is called a singular Sturm-Liouville problem \cite{Arfken}. Sturm-Liouville problems appear in various branches of engineering, physics and biology. Recently in \cite{theBook}, some generalized Sturm-Liouville problems in three different continuous, discrete and q-discrete spaces have been introduced and classified.

Let $ y_{n}(x) $ and $ y_{m}(x) $ be two solutions of equation \eqref{EQ1.ODE}. Following Sturm-Liouville theory \cite{Arfken, Niki}, these functions are orthogonal with respect to the positive weight function $ \rho(x) $ on $ (a,b) $ under the given conditions \eqref{EQ2.boundary}, i.e.
\begin{equation}\label{EQ3.defOrtho}
\int_a^b {\rho (x){y_n}(x){y_m}(x)\,dx}  = \left( {\int_a^b {\rho (x)y_n^2(x)\,dx} } \right){\delta _{n,m}},
\end{equation}
where
\[{\delta _{n,m}} = \left\{ \begin{array}{l}
0\,\,\,\,\,\,\,\,(n \ne m),\\
1\,\,\,\,\,\,\,\,(n = m).
\end{array} \right.\]

Many special functions in theoretical and mathematical physics are solutions of a regular or singular Sturm-Liouville problem satisfying the orthogonality condition \eqref{EQ3.defOrtho}, see e.g. \cite{Arfken, MR0481884, MMS, Szego}.

There are totally six sequences of real polynomials \cite{MR2246501,MR1915513} that are orthogonal with respect to the Pearson distributions family
\begin{equation}\label{eq:3}
W\left(
\begin{array}{r|l}
\begin{array}{cc}
{d^{*},\,\,e^{*}}\\
{a,\,\,b,\,\,c}
\end{array} & 
{x}
\end{array} \right) = \exp \left( \int \frac{d^{*}x+e^{*}}{ax^2+bx+c}\,dx \right) \qquad (a,b,c,d^{*},e^{*} \in {\mathbb{R}}).
\end{equation}
Three of them (i.e. Jacobi, Laguerre and Hermite polynomials \cite{MR1149380}) are infinitely orthogonal with respect to three special cases of the positive function \eqref{eq:3} (i.e. Beta, Gamma and Normal distributions \cite{Kotz}) and three other ones are finitely orthogonal with respect to F-Fisher, Inverse Gamma and Generalized T-Student distributions \cite{MR2053407,Kotz}, limited to some parametric constraints. Table \ref{table1} shows the main properties of these six sequences.

\begin{table}[h]
\caption{Characteristics of six sequences of orthogonal polynomials} \label{table1} \centering %
\begin{tabular}{|c|c|c|c|} 
\hline 
\begin{tabular}{c} Polynomial \\ notation \end{tabular} & Distribution & Weight function & \begin{tabular}{c} Kind, \\ Interval and \\ Parameters constraint \end{tabular}  \\ \hline
 $P_n^{(\alpha,\beta)} (x)$ & Beta & 
\begin{tabular}{c} 
$W\left(\begin{array}{r|l}\begin{array}{c}-\alpha-\beta,\,\,-\alpha+\beta\\-1,\,\,0,\,\,1\end{array} & {x}\end{array} \right)$ \\[4mm]
$=(1 - x)^\alpha (1 + x)^\beta $ \end{tabular} &  \begin{tabular}{c} Infinite \\ $[-1,1]$ \\ $\forall n, \alpha>-1, \beta>-1$ \end{tabular}  \\ \hline
 $L_n^{(\alpha)} (x)$ & Gamma & 
\begin{tabular}{c}
$W\left(\begin{array}{r|l}\begin{array}{c}-1,\,\,\alpha\\0,\,\,1,\,\,0\end{array} & {x}\end{array} \right)$ \\[4mm]
$=x^\alpha \exp ( - x)$ 
\end{tabular}
& \begin{tabular}{c} Infinite \\ $[0,\infty)$ \\ $\forall n, \alpha>-1$ \end{tabular} \\ \hline
 $H_n (x)$ & Normal & 
\begin{tabular}{c} 
$W\left(\begin{array}{r|l}\begin{array}{c}-2,\,\,0\\0,\,\,0,\,\,1\end{array} & {x}\end{array} \right)$ \\[4mm]
$=\exp ( - x^2 )$ 
\end{tabular}
& \begin{tabular}{c} Infinite \\ $(-\infty,\infty)$ \\  --- \end{tabular}    \\ \hline
$M_n^{(p,q)} (x)$ & Fisher F  & 
\begin{tabular}{c} 
$W\left(\begin{array}{r|l}\begin{array}{c}-p,\,\,q\\1,\,\,1,\,\,0\end{array} & {x}\end{array} \right)$ \\[4mm]
$=x^q (x + 1)^{ - (p + q)}$ 
\end{tabular}
& \begin{tabular}{c} Finite \\ $[0,\infty)$ \\ $\max n<(p-1)/2$ \\ $q>-1$ \end{tabular}  \\ \hline
$N_n^{(p)} (x)$ & Inverse Gamma  & 
\begin{tabular}{c} 
$W\left(\begin{array}{r|l}\begin{array}{c}-p,\,\,1\\1,\,\,0,\,\,0\end{array} & {x}\end{array} \right)$ \\[4mm]
$=x^{ - p} \exp ( - 1/x)$ 
\end{tabular}
& \begin{tabular}{c} Finite \\ $[0,\infty)$ \\ $\max n<(p-1)/2$ \end{tabular}  \\ \hline
 $J_n^{(p,q)} (x)$ & Generalized T  & 
\begin{tabular}{c} $W\left(\begin{array}{r|l}\begin{array}{c}-2p,\,\,q\\1,\,\,0,\,\,1\end{array} & {x}\end{array} \right)$ \\[4mm] 
$=(1+x^2)^{-p} \exp(q \arctan x)$ \end{tabular} & \begin{tabular}{c} Finite \\ $(-\infty,\infty)$ \\ $\max n <p-1/2$ \end{tabular}\\ \hline 
\end{tabular}
\end{table}

It has been proved by S. Bochner \cite{MR1545034,MR0481884} that if an infinite sequence of polynomials $ \lbrace P_{n}\rbrace_{n=0}^{\infty} $ satisfies a second-order eigenvalue equation of the form
\begin{equation*}
\sigma (x){P''_n}(x) + \tau (x){P'_n}(x)+r(x)P_{n}(x) = {\lambda _n}{P_n}(x) \qquad \qquad n=0,1,2,\ldots,
\end{equation*}
then $ \sigma (x) $, $ \tau (x) $ and $ r(x) $ must be polynomials of degree $ 2,\,1 $ and $ 0 $, respectively. Moreover, if the sequence $ \lbrace P_{n}\rbrace_{n=0}^{\infty} $ is orthogonal, then it has to be one of the classical Jacobi, Laguerre or Hermite polynomials which satisfy a second order differential equation of the form \cite{AlSalam, Andrews, MR1545034}
\begin{equation}\label{EQ.hypDiff}
\sigma (x){y''_n}(x) + \tau (x){y'_n}(x) - {\lambda _n}{y_n}(x) = 0{\rm{ ,}}
\end{equation}
where 
\[ \sigma (x) = a{x^2} + bx + c \,\,\,\,\text{ and }\,\,\,\, \tau (x) = dx + e, \] 
and 
\[ {\lambda _n} = n\big(d+(n-1)a\big),\]
is the eigenvalue depending on $ n=0,1,2,\ldots $. However, there are three other sequences of hypergeometric polynomials that are solutions of the equation \eqref{EQ.hypDiff} but finitely orthogonal \cite{MR2053407, MR1915513}.

\medskip

It is usually supposed in the literature that the orthogonal polynomial systems start with a polynomial of degree 0. Nevertheless, from Sturm-Liouville theory point of view, this restriction is not necessary \cite{MR2542180}.

Recently in \cite{MR2542180,MR2610341}, two new families of exceptional orthogonal polynomials, i.e. Jacobi X$_{1}$-polynomials $\hat{P}_{n}^{(\alpha,\beta)}(x)$ and Laguerre X$_{1}$-polynomials $\hat{L}_{n}^{(\alpha)}(x)$ have been introduced as solutions of a second-order eigenvalue equation of the form
\begin{multline}
 \Big( k_{2}(x-b)^{2} + k_{1}(x-b) + k_{0} \Big) y''_{n}(x) +  \frac{ax-ab-1}{x-b} \Big( k_{1}(x-b) + 2 k_{0} \Big) y'_{n}(x) \\
-\Big(  \frac{a}{x-b}\big( k_{1}(x-b)+2k_{0} \big) + \lambda_{n} \Big) y_{n}(x)=0,
\end{multline}
for $n \geq 1$, where 
\[
\lambda_{n}=(n-1)(n k_{2}-ak_{1}),
\]
and $ k_{0}\neq 0, k_{1} $, $ k_{2} $ are real constants. In this sense, the authors in \cite{MR2542180} have proved a converse statement similar to Bochner's theorem for the classical orthogonal polynomials: if a self-adjoint second order operator has a polynomial eigenfunctions $\{P_{i}(x)\}_{i=1}^{\infty}$, then it must be either the X$_{1}$-Jacobi or the X$_{1}$-Laguerre Sturm-Liouville problem.
Moreover, the functions
\begin{equation}\label{eq:Wxjacobi}
\hat{W}_{\alpha,\beta}(x)={\Big(x-\dfrac{\beta+\alpha}{\beta-\alpha}\Big)^{-2}}{(1-x)^{\alpha}(1+x)^{\beta}}\qquad\text{for} \quad x\in(-1,1),
\end{equation}
with the restrictions $ \alpha,\beta>-1,\,\alpha\neq\beta,\,\text{sgn}\,\alpha=\text{sgn}\,\beta $,
and 
\begin{equation}\label{eq:Wxlaguerre}
\hat{W}_{\alpha}(x)=\big(x+\alpha\big)^{-2}x^{\alpha}e^{-x}\qquad\text{for} \quad x\in(0,\infty),
\end{equation}
with the restriction $ \alpha>0 $, are the weight functions corresponding to $\hat{P}_{n}^{(\alpha,\beta)}(x)$ and $\hat{L}_{n}^{(\alpha)}(x)$, respectively.

Exceptional orthogonal polynomials have been recently of great interest due to their important applications in exactly solvable potentials and supersymmetry \cite{MR2439200}, Dirac operators minimally coupled to external fields \cite{MR2771725} and entropy measures in quantum information theory \cite{MR2905628}. Moreover, the relationship between exceptional orthogonal polynomials and Darboux transformations has been observed, giving rise to new families of X$_{2}$ polynomials of codimension two \cite{MR2439200,MR2559677}. See also \cite{MR2569488,MR2588057} for higher-order codimensional families.

In this paper, we consider six sequences of X$ _{1} $ orthogonal polynomials as special solutions of a generic Sturm-Liouville equation of the form
\begin{multline}\label{eq:5}
(x-r)\left( a_{2} x^{2} + a_{1} x + a_{0} \right) y_{n}''(x) + \left( b_{2} x^{2}+ b_{1} x+ b_{0} \right) y_{n}'(x) \\
- \big( \lambda_{n}(x-r)+c_{0}^{*}\big) y_{n}(x)=0,\qquad  n \geq 1,
\end{multline}
where $ r $ is a real parameter such that $a_{2} r^{2} + a_{1} r + a_{0} \neq 0 $ and the roots of $ b_{2} x^{2}+ b_{1} x+ b_{0} $ are supposed to be real.

Both infinite and finite types of exceptional orthogonal X$ _{1} $-polynomials can be extracted from the above equation \eqref{eq:5}. 
Although three infinite polynomial sequences have been investigated in \cite{Ferrero} for only some particular parameters, the finite cases of exceptional X$ _{1} $-polynomials orthogonal with respect to three particular weight functions on infinite intervals are introduced in this paper for the first time. 
A fundamental point is that the weight functions corresponding to these six sequences are exactly a multiplication of Pearson distributions family introduced in table \ref{table1}.

\section{A Review on Classical Orthogonal Polynomials}

It is shown in \cite{MR2246501} that the monic polynomial solution of equation \eqref{EQ.hypDiff} can be represented as  
\begin{equation}\label{Dissertation.EQ.6.8}
{y_n}(x) = {\bar{P}_n}\left( {\left. {\begin{array}{*{20}{c}}
{\begin{array}{*{20}{c}}
d,&e
\end{array}}\\
{a\,,\,\,b\,,\,\,c}
\end{array}} \right|x} \right) = \sum\limits_{k = 0}^n {\,\binom{n}{k}G_k^{(n)}(a,b,c,d,e)} \,{x^k},
\end{equation}
where 
\[G_k^{(n)} = {\left( {\frac{{2a}}{{b + \sqrt {{b^2} - 4ac} }}} \right)^{k - n}}{}_2{F_1}\left( {\begin{array}{*{20}{c}}
{k- n,\begin{array}{*{20}{c}}
{}&{\frac{{2ae - bd}}{{2a\sqrt {{b^2} - 4ac} }} + 1 - \frac{d}{{2a}} - n}
\end{array}}\\[3mm]
{-\frac{d}{a}+2- 2n}
\end{array}\left| {\,\frac{{2\sqrt {{b^2} - 4ac} }}{{b + \sqrt {{b^2} - 4ac} }}} \right.} \right),\]
and
\[ _2{F_1}\left( {\left.  \begin{array}{*{20}{c}}
{\begin{array}{*{20}{c}}
{a,}&{b}
\end{array}}\\
{c}
\end{array}\right| {x} } \right)=\sum\limits_{k=0}^{\infty}\dfrac{(a)_{k}(b)_{k}}{(c)_{k}}\dfrac{x^{k}}{k!} ,\] 
is the Gauss hypergeometric function \cite{Slater} for $ {(a)_k} = a(a+1)\ldots (a+k-1) $.

The general formula \eqref{Dissertation.EQ.6.8} is a suitable tool to compute the coefficients of $ x^k $ for any fixed degree $ k $ and arbitrary $ a $. For example, to obtain the coefficient $ x^{n-1} $, it is enough to calculate the term
\begin{multline}\label{Dissertation.EQ.6.21}
G_{n - 1}^{(n)}(a,b,c,d,e) = {(\frac{{2a}}{{b + \Delta }})^{ - 1}}\,{}_2{F_1}\left( {\left. \begin{array}{*{20}{c}}
{\begin{array}{*{20}{c}}
{ - 1,}&{\frac{{2ae - bd - (d + (2n - 2)a)\Delta }}{{2a\Delta }}}
\end{array}}\\[3mm]
{2-2n-\frac{d}{a}}
\end{array}\right|\,\frac{{2\Delta }}{{b + \Delta }}} \right)\\
 = \big(\frac{{b + \Delta }}{{2a}}\big)\Big(1 + \frac{{2ae - bd - (d + (2n - 2)a)\Delta }}{{2a\Delta }}\, \frac{{2\Delta }}{{b + \Delta }} \, \frac{a}{{d + (2n - 2)a}}\Big)\\
 = \frac{{e + (n - 1)b}}{{d + (2n - 2)a}},
\end{multline}
in which $ \Delta  = \sqrt {{b^2} - 4ac} $. Note in \eqref{Dissertation.EQ.6.21} that all parameters are free and can adopt any value including zero because neither both values $ a $ and $ d $ nor both values $ b $ and $ e $ can vanish together in \eqref{Dissertation.EQ.6.8}. After simplifying $ G_k^{(n)}(a,b,c,d,e) $ for $ k=n-1,n-2,\ldots $ we eventually obtain
\begin{multline*}%\label{Dissertation.EQ.6.22}
{{\bar P}_n}\left( {\left. \begin{array}{*{20}{c}}
{\begin{array}{*{20}{c}}
d&e
\end{array}}\\
{a\,\,\,b\,\,\,c}
\end{array}\right| {\,x} } \right) = {x^n} + \left( {\begin{array}{*{20}{c}}
n\\
1
\end{array}} \right)\frac{{e + (n - 1)b}}{{d + (2n - 2)a}}{x^{n - 1}} \\[3mm]
+ \left( {\begin{array}{*{20}{c}}
n\\
2
\end{array}} \right)\frac{{(e + (n - 1)b)(e + (n - 2)b) + c(d + (2n - 2)a)}}{{(d + (2n - 2)a)(d + (2n - 3)a)}}{x^{n - 2}}
 + ... \\[3mm]
 + \left( {\begin{array}{*{20}{c}}
n\\
n
\end{array}} \right)\,{\Big(\frac{{b + \sqrt {{b^2} - 4ac} }}{{2a}}\Big)^n}{}_2{F_1}\left( {\left. \begin{array}{*{20}{c}}
{\begin{array}{*{20}{c}}
{ - n,}&{\frac{{2ae - bd - (d + (2n - 2)a)\sqrt {{b^2} - 4ac} }}{{2a\sqrt {{b^2} - 4ac} }}}
\end{array}}\\[3mm]
{2-2n-\frac{d}{a}}
\end{array}\right|\,\frac{{2\sqrt {{b^2} - 4ac} }}{{b + \sqrt {{b^2} - 4ac} }}} \right){\rm{ }}.
\end{multline*}
For instance, we have
\begin{align*}
{{\bar P}_0}\left( {\left. \begin{array}{*{20}{c}}
{\begin{array}{*{20}{c}}
d&e
\end{array}}\\
{a\,\,\,b\,\,\,c}
\end{array}\right| {\,x}} \right) &= 1\,\,,\nonumber\\
{{\bar P}_1}\left( {\left. \begin{array}{*{20}{c}}
{\begin{array}{*{20}{c}}
d&e
\end{array}}\\
{a\,\,\,b\,\,\,c}
\end{array}\right| {\,x} } \right) &= x + \frac{e}{d}\,\,,\nonumber\\
{{\bar P}_2}\left( {\begin{array}{*{20}{c}}
{\begin{array}{*{20}{c}}
d&e
\end{array}}\\
{a\,\,\,b\,\,\,c}
\end{array}\left| {\,x} \right.} \right) &= {x^2} + 2\frac{{e + b}}{{d + 2a}}x + \frac{{c(d + 2a) + e(e + b)}}{{(d + 2a)(d + a)}}\,\,,\nonumber\\
{{\bar P}_3}\left( {\left. \begin{array}{*{20}{c}}
{\begin{array}{*{20}{c}}
d&e
\end{array}}\\
{a\,\,\,b\,\,\,c}
\end{array}\right| {\,x} } \right) &= {x^3} + 3\frac{{e + 2b}}{{d + 4a}}{x^2} + 3\frac{{c(d + 4a) + (e + b)(e + 2b)}}{{(d + 4a)(d + 3a)}}x\nonumber\\
& \quad + \frac{{2c(d + 3a)(e + 2b) + ce(d + 4a) + e(e + b)(e + 2b)}}{{(d + 4a)(d + 3a)(d + 2a)}}\,.%\label{Dissertation.EQ.6.22.2}
\end{align*}

Moreover, by referring to the Nikiforov and Uvarov approach \cite{Niki} and considering equation \eqref{EQ.hypDiff} as a self-adjoint form, the Rodrigues representation of the monic polynomials is derived as
\begin{multline}\label{Dissertation.EQ.6.24}
{\bar P_n}\left( {\left. \begin{array}{*{20}{c}}
{\begin{array}{*{20}{c}}
d&e
\end{array}}\\
{a\,\,\,b\,\,\,c}
\end{array}\right| {\,x} } \right) = \frac{1}{{\big(\prod\limits_{k = 1}^n {d + (n + k - 2)a} \big)\,W \left( {\left. \begin{array}{*{20}{c}}
{\begin{array}{*{20}{c}}
d&e
\end{array}}\\
{a\,\,\,b\,\,\,c}
\end{array}\right| {\,x}} \right)}}\\
 \times \frac{{{d^n}\Big({{(a{x^2} + bx + c)}^n}W \left( {\left. \begin{array}{*{20}{c}}
{\begin{array}{*{20}{c}}
d&e
\end{array}}\\
{a\,\,\,b\,\,\,c}
\end{array}\right| {\,x} } \right)\Big)}}{{d{x^n}}}\,,
\end{multline}
where 
\begin{equation}\label{EQ.3scan1.4.22}
W \left( {\left. \begin{array}{*{20}{c}}
{\begin{array}{*{20}{c}}
d&e
\end{array}}\\
{a\,\,\,b\,\,\,c}
\end{array}\right| {\,x} } \right) = \exp \Big(\int {\frac{{(d - 2a)x + e - b}}{{a{x^2} + bx + c}}dx} \Big) .
\end{equation}

By using the formulas \eqref{Dissertation.EQ.6.8} or \eqref{Dissertation.EQ.6.24} we can also obtain a generic three term recurrence equation for the polynomials as follows \cite{MR2246501}
\begin{multline*}%\label{Dissertation.EQ.6.30}
{{\bar P}_{n + 1}}(x) = \left( {x + \frac{{2n(n + 1)ab + (d - 2a)(e + 2nb)}}{{(d + 2na)(d + (2n - 2)a)}}} \right)\,{{\bar P}_n}(x)\\
 +n(d + (n - 2)a)\times\\
  \frac{{\left( {c{{(d + (2n - 2)a)}^2} - n{b^2}(d + (n - 2)a) + (e - b)(a(e + b) - bd)} \right)}}{{(d + (2n - 3)a){{(d + (2n - 2)a)}^2}(d + (2n - 1)a)}}{{\bar P}_{n - 1}}(x),
\end{multline*}
in which $ {{\bar P}_{n}}(x) $ denotes the same as monic polynomials of \eqref{Dissertation.EQ.6.8} with the initial values 
\[ {\rm{ }}{\bar P_0}(x) = 1\,\,\,\,{\rm{and}}\,\,\,\,{\bar P_1}(x) = x + \frac{e}{d} .\]

\medskip

Finally, the norm square value of the monic polynomials \eqref{Dissertation.EQ.6.8} can be simplified and computed as follows:\\
Let $ [L,U] $ be a predetermined orthogonality interval which consists of the zeros of $ \sigma (x) = a{x^2} + bx + c $ or $ \pm\infty $. By noting the Rodrigues representation \eqref{Dissertation.EQ.6.24} we have
\begin{multline}\label{Dissertation.EQ.6.31}
{\left\| {{{\bar P}_n}} \right\|^2} = \int_L^U {\,{{\bar P}_n}^2\left( {\left. \begin{array}{*{20}{c}}
{\begin{array}{*{20}{c}}
d&e
\end{array}}\\
{a\,\,\,b\,\,\,c}
\end{array}\right| {\,x} } \right)\,W\left( {\left. \begin{array}{*{20}{c}}
{\begin{array}{*{20}{c}}
d&e
\end{array}}\\
{a\,\,\,b\,\,\,c}
\end{array}\right| {\,x}} \right)dx}  = \frac{1}{{\prod\limits_{k = 1}^n {d + (n + k - 2)a} }}\\
\times \int_L^U {\,{{\bar P}_n}\left( {\left. \begin{array}{*{20}{c}}
{\begin{array}{*{20}{c}}
d&e
\end{array}}\\
{a\,\,\,b\,\,\,c}
\end{array}\right| {\,x} } \right)\,\dfrac{d^n}{d{x^n}}\Big({{(a{x^2} + bx + c)}^n}W\left( {\left. \begin{array}{*{20}{c}}
{\begin{array}{*{20}{c}}
d&e
\end{array}}\\
{a\,\,\,b\,\,\,c}
\end{array}\right| {\,x} } \right)\Big)\,dx} .
\end{multline}
So, integrating by parts from the right hand side of \eqref{Dissertation.EQ.6.31} yields 
\begin{equation*}\label{Dissertation.EQ.6.32}
{\left\| {{{\bar P}_n}} \right\|^2} = \frac{{n!\,{{( - 1)}^n}}}{{\prod\limits_{k = 1}^n {d + (n + k - 2)a} }}\int_L^U {{{(a{x^2} + bx + c)}^n}\Big(\exp \int {\frac{{(d - 2a)x + e - b}}{{a{x^2} + bx + c}}dx}\Big )\,dx} {\rm{ }}{\rm{.}}
\end{equation*}

\medskip

Although the Jacobi polynomials \[ \bar{P}^{(\alpha,\beta)}_{n}(x)=
{\bar{P}_n}\left( {\left. {\begin{array}{*{20}{c}}
{\begin{array}{*{20}{c}}
-\alpha-\beta-2, & \beta-\alpha
\end{array}}\\
{-1,\,\,\,0,\,\,\,1}
\end{array}} \right|x} \right), \]
Laguerre polynomials
\[\bar{L}^{(\alpha)}_{n}(x)=
{\bar{P}_n}\left( {\left. {\begin{array}{*{20}{c}}
{\begin{array}{*{20}{c}}
-1, & \alpha+1
\end{array}}\\
{0,\,\,\,1,\,\,\,0}
\end{array}} \right|x} \right),  \]
and Hermit polynomials 
\[\bar{H}_{n}(x)={\bar{P}_n}\left( {\left. {\begin{array}{*{20}{c}}
{\begin{array}{*{20}{c}}
-2 ,& 0
\end{array}}\\
{0,\,\,\,0,\,\,\,1}
\end{array}} \right|x} \right), \]
are three polynomial solutions of equation \eqref{EQ.hypDiff}, there are three other sequences of hypergeometric polynomials that are finitely orthogonal with respect to the generalized T, inverse Gamma and F distributions \cite{MR2053407, MR1915513} and are solutions of equation \eqref{EQ.hypDiff}.

The first finite sequence of classical orthogonal polynomials as 
\begin{equation}\label{EQ.defM}
\bar{M}_n^{(p,q)}(x) ={\bar{P}_n}\left( {\left. {\begin{array}{*{20}{c}}
{\begin{array}{*{20}{c}}
2-p , & 1+q
\end{array}}\\
{1,\,\,\,1,\,\,\,0}
\end{array}} \right|x} \right),
\end{equation}
satisfies the differential equation
\begin{equation}\label{EQ.diF1}
({x^2} + x){y''_n}(x) + \big((2 - p)x + q+1\big){y'_n}(x) - n(n + 1 - p){y_n}(x) = 0,
\end{equation}
and is finitely orthogonal with respect to the weight function 
\[ {W_1}(x;p,q) = {x^q}{(1 + x)^{ - (p + q)}} ,\] 
on $ [0,\infty) $ if and only if \cite{MR1915513}
\[p > 2\{ \max \,n\} \, + 1 \quad\text{and} \quad q>-1.\]

The second finite sequence defined as
\begin{equation*}
\bar{N}_n^{(p)}(x)=
{\bar{P}_n}\left( {\left. {\begin{array}{*{20}{c}}
{\begin{array}{*{20}{c}}
2-p, & 1
\end{array}}\\
{1,\,\,\,0,\,\,\,0}
\end{array}} \right|x} \right),
\end{equation*}
satisfies the differential equation
\begin{equation}\label{EQ.diF2}
{x^2}{y''_n}(x) + \big((2 - p)x + 1\big){y'_n}(x) - n(n + 1 - p){y_n}(x) = 0,
\end{equation}
and is finitely orthogonal with respect to the weight function \cite{MR1915513} 
\[W_{2}(x;p) =  {x^{- p}}{e^{ - \frac{1}{x}}},\] 
on $ (0,\infty) $ for $ n = 0,1,2,...,N < \frac{{p - 1}}{2} $.

Finally, the third finite sequence, which is finitely orthogonal with respect to the generalized T-student distribution weight function
\[ W_{3}(x;p,q) = {\left( 1+x^2 \right)^{ - p}}\exp (q\arctan x),\] 
is defined on $ (-\infty,\infty) $ as 
\begin{equation*}
\bar{J}_n^{(p,q)}(x)=
{\bar{P}_n}\left( {\left. {\begin{array}{*{20}{c}}
{\begin{array}{*{20}{c}}
2-2p, & q
\end{array}}\\
{1,\,\,\,0,\,\,\,1}
\end{array}} \right|x} \right),
\end{equation*}
satisfying the equation 
\begin{equation}\label{EQ.diF32}
(1+{x^2} )\,y''_{n}(x) + \big( 2(1-p)x+q \big)\,{y'_n}(x)
 - n(n + 1 - 2p)\,{y_n}(x) = 0,
\end{equation}
and the orthogonality property holds if
\[ n = 0,1,2,...,N < p - \frac{1}{2}\,\,\,\text{and}\,\,q \in {\mathbb{R}} .\]

\section{A Generic Classification of Exceptional X$ _{1} $-Orthogonal Polynomials}

Since for any arbitrary real parameters $ A,B $ and $ C $ the relation
\begin{equation*}
Ax^{2}+Bx+C=A(x-r)^{2}+(2Ar+B)(x-r)+Ar^{2}+Br+C,
\end{equation*}
always holds true, another form of equation \eqref{eq:5} is as 
\begin{multline}\label{eq:6}
(x-r)\Big(a_{2}(x-r)^{2} + (2a_{2}r+a_{1})(x-r) + a_{2}r^{2}+a_{1}r+a_{0} \Big) y_{n}''(x)\\ + \Big( b_{2} (x-r)^{2} + (2b_{2}r+b_{1})(x-r) + b_{2}r^{2}+b_{1}r+b_{0} \Big) y_{n}'(x) \\
-\Big(\lambda_{n}(x-r)+c_{0}^{*}  \Big)y_{n}(x)=0,\qquad  n \geq 1.
\end{multline}

The eigenvalue $\lambda_{n}$ in \eqref{eq:6} is to be determined such that for every $ n\geq 1 $, the solution $ y_{n} $ is a polynomial of degree $ n $. 
For this purpose, we first consider a subspace of the whole space of polynomials of degree at most $ n $ as 
\[\Pi_{n,r,\nu}=\text{span}\Big\lbrace (x-r-{\nu}),(x-r)^{2},\ldots, (x-r)^{n}\Big\rbrace,\]
in which $ \nu $ is a real constant. Then substituting $ y_{1}(x)=x-r-{\nu} $ and $ y_{n}(x)=(x-r)^{n} $ for $ n\geq 2 $ into \eqref{eq:6} respectively yield
\begin{equation}\label{eq:bc}
\big(b_{2}-\lambda_{1}\big)(x-r)^{2}+\big(2b_{2}r+b_{1}-c_{0}^{*}+{\nu}\lambda_{1}\big)(x-r)+b_{2}r^{2}+b_{1}r+b_{0}+{\nu}c_{0}^{*}=0,
\end{equation}
and
\begin{multline*}
n(n-1)\Big(a_{2}(x-r)^{2} + (2a_{2}r+a_{1})(x-r) + a_{2}r^{2}+a_{1}r+a_{0} \Big)(x-r)^{n-1}\\
+n\Big( b_{2} (x-r)^{2} + (2b_{2}r+b_{1})(x-r) + b_{2}r^{2}+b_{1}r+b_{0} \Big)(x-r)^{n-1}\\
-\Big(\lambda_{n}(x-r)+c_{0}^{*}\Big)(x-r)^{n}=0\qquad\quad n\geq 2.
\end{multline*}
Therefore
\begin{equation*}
\lambda_{n}=n\big((n-1)a_{2}+b_{2}\big)\qquad \text{for}\qquad n\geq 1,
\end{equation*}
and
\begin{equation}\label{eq:2.2.1}
\begin{cases}
2b_{2}r+b_{1}-c_{0}^{*}+\nu b_{2}=0,\\[2mm]
b_{2}r^{2}+b_{1}r+b_{0}+\nu c_{0}^{*}=0.
\end{cases}
\end{equation}
By solving the system \eqref{eq:2.2.1} we get
\begin{equation*}
\nu=\dfrac{-(2b_{2}r+b_{1})\pm\sqrt{b_{1}^{2}-4b_{0}b_{2}}}{2b_{2}}=\begin{cases}
r_{1}-r,\\
r_{2}-r,\end{cases} 
\end{equation*}
where $ r_{1},r_{2}$ are roots of $ b_{2}x^{2}+b_{1}x+b_{0} $, and
\[c_{0}^{*}=\dfrac{2b_{2}(b_{2}r^{2}+b_{1}r+b_{0})}{2b_{2}r+b_{1}\mp\sqrt{b_{1}^{2}-4b_{0}b_{2}}}=
\begin{cases}
b_{2}(r-r_{2}),\\[4mm]
b_{2}(r-r_{1}).
\end{cases}\]

\begin{corollary}\label{coro1}
If we take $ b_{2}x^{2}+b_{1}x+b_{0}=b_{2}(x-r_{1})(x-r_{2}) $ and 
\[ \Pi_{n,r,\nu}=\text{span}\lbrace e_{k}(x)\rbrace_{k=1}^{n} ,\] 
then\\
(i) $ e_{1}(x)=x-r_{1} $ and $ \big\lbrace e_{k}(x)\big\rbrace_{k=2}^{\infty}=\big\lbrace (x-r)^{k}\big\rbrace_{k=2}^{\infty} $ lead to $ c_{0}^{*}=b_{2}(r-r_{2}) $.\\
\\
(ii) $ e_{1}(x)=x-r_{2} $ and $ \big\lbrace e_{k}(x)\big\rbrace_{k=2}^{\infty}=\big\lbrace (x-r)^{k}\big\rbrace_{k=2}^{\infty} $ lead to $ c_{0}^{*}=b_{2}(r-r_{1}) $.\\
\\
We will use this important corollary in the next sections. 
\end{corollary}

We are now in a good position to prove that the polynomial solutions of equation \eqref{eq:6} in $ \Pi_{n,r,\nu} $ are orthogonal on an interval, say $ [a,b] $, with respect to a weight function in the form
\begin{equation}\label{eq:7}
\rho(x)=(x-r) \omega(x),
\end{equation}
where $\omega(x)$ satisfies the differential equation
\begin{equation}\label{eq:8}
\frac{\omega'(x)}{\omega(x)} = 
\frac{(b_{2}-3a_{2})x^{2}+(b_{1}-2a_{1}+2a_{2}r)x+b_{0}-a_{0}+a_{1}r}{(x-r) \big(a_{2}x^{2}+a_{1}x+a_{0} \big)}.
\end{equation}
To prove this claim, we first consider the self-adjoint form of equation \eqref{eq:6} as \begin{equation}\label{eq:sn}
\Big(\omega(x)(x-r)\big(a_{2}x^{2} + a_{1}x + a_{0} \big)y'_{n}\Big)'=\omega(x)\Big(\lambda_{n}(x-r)+c_{0}^{*}  \Big)y_{n}(x),
\end{equation}
and for the index $ m $ as
\begin{equation}\label{eq:sm}
\Big(\omega(x)(x-r)\big(a_{2}x^{2} + a_{1}x + a_{0} \big)y'_{m}\Big)'=\omega(x)\Big(\lambda_{m}(x-r)+c_{0}^{*}  \Big)y_{m}(x).
\end{equation}
Multiplying by $ y_{m}(x) $ and $  y_{n}(x) $ in relations \eqref{eq:sn} and \eqref{eq:sm} respectively, subtracting them and then integrating from both sides we get
\begin{multline}\label{eq:ortho1}
\Big[ \omega(x)(x-r)\big(a_{2}x^{2} + a_{1}x + a_{0} \big)\big(y'_{n}(x)y_{m}(x)-y'_{m}(x)y_{n}(x)\big)\Big]_{a}^{b}\\
=(\lambda_{n}-\lambda_{m})\int_{a}^{b}(x-r)\omega(x)y_{n}(x)y_{m}(x)\,dx.
\end{multline}
Now if the following relations 
\begin{align*}
\omega(a)(a-r)\big(a_{2}a^{2} + a_{1}a + a_{0} \big)&=0,\\
\omega(b)(b-r)\big(a_{2}b^{2} + a_{1}b + a_{0} \big)&=0,
\end{align*}
hold, the left hand side of \eqref{eq:ortho1} is equal to zero and therefore
\begin{equation*}
\int_{a}^{b}(x-r)\omega(x)y_{n}(x)y_{m}(x)\,dx=0\qquad m\neq n,
\end{equation*}
which shows the orthogonality of polynomial sequence $ \left\lbrace y_{n}(x) \right\rbrace_{n=1}^{\infty}  $ with respect to the weight function $ \rho(x)=(x-r)\omega(x) $. On the other hand, according to \eqref{eq:8}, 
$ \rho(x) $ should have ``seven'' free parameters ( regardless $ r $, one parameter more than the number of parameters in Pearson distribution) because the explicit solution of equation \eqref{eq:8} is as
\begin{equation}\label{eq:9}
\omega(x)=\exp \left( \int\frac{(b_{2}-3a_{2})x^{2}+(b_{1}-2a_{1}+2a_{2}r)x+b_{0}-a_{0}+a_{1}r}{(x-r) \big(a_{2}x^{2}+a_{1}x+a_{0} \big)}\,dx \right)\,.
\end{equation}

Here the key point is that the function \eqref{eq:9} is exactly a multiplication of Pearson distribution given in \eqref{eq:3}, because if the integrand function of \eqref{eq:9} is written as a sum of two fractions with linear and quadratic denominators in the form
\begin{multline*}
\dfrac{(b_{2}-3a_{2})x^{2}+(b_{1}-2a_{1}+2a_{2})x+b_{0}-a_{0}+a_{1}r}{(x-r) \left( a_{2}x^{2}+a_{1}x+a_{0} \right)}  
= \dfrac{\,\,\frac{b_{2}r^{2}+b_{1}r+b_{0}}{a_{2}r^{2}+a_{1}r+a_{0}}-1\,\,}{x-r} \\[3mm]
+ \dfrac{\big(b_{2}-a_{2}(2+\frac{b_{2}r^{2}+b_{1}r+b_{0}}{a_{2}r^{2}+a_{1}r+a_{0}})\big)x+b_{1}+b_{2}r-(\frac{b_{2}r^{2}+b_{1}r+b_{0}}{a_{2}r^{2}+a_{1}r+a_{0}})(a_{1}+a_{2}r)-a_{1}}{a_{2}x^{2}+a_{1}x+a_{0}},
\end{multline*}
then we obtain
\begin{align}
&\omega(x)=(x-r)^{\frac{b_{2}r^{2}+b_{1}r+b_{0}}{a_{2}r^{2}+a_{1}r+a_{0}}-1}\\ 
&\quad\times \exp \left( \int\dfrac{\big(b_{2}-a_{2}(2+\frac{b_{2}r^{2}+b_{1}r+b_{0}}{a_{2}r^{2}+a_{1}r+a_{0}})\big)x+b_{1}+b_{2}r-(\frac{b_{2}r^{2}+b_{1}r+b_{0}}{a_{2}r^{2}+a_{1}r+a_{0}})(a_{1}+a_{2}r)-a_{1}}{a_{2}x^{2}+a_{1}x+a_{0}}\,dx \right) \nonumber\\[3mm]
&\,\,\,\quad=(x-r)^{\frac{b_{2}r^{2}+b_{1}r+b_{0}}{a_{2}r^{2}+a_{1}r+a_{0}}-1}\nonumber\\
&\quad\times \,W\left(
\begin{array}{r|l}
\begin{array}{cc}
{b_{2}-a_{2}(2+\frac{b_{2}r^{2}+b_{1}r+b_{0}}{a_{2}r^{2}+a_{1}r+a_{0}}),\,\,b_{1}+b_{2}r-(\frac{b_{2}r^{2}+b_{1}r+b_{0}}{a_{2}r^{2}+a_{1}r+a_{0}})(a_{1}+a_{2}r)-a_{1}}\\[2mm]
{a_{2}\,,\,a_{1}\,,\,\,a_{0}}
\end{array} & 
{x}
\end{array} \right),
\nonumber
\end{align}
and accordingly,
\begin{equation}\label{eq:10}
\rho(x)=(x-r)^{\frac{b_{2}r^{2}+b_{1}r+b_{0}}{a_{2}r^{2}+a_{1}r+a_{0}}}\,W\left(
\begin{array}{r|l}
\begin{array}{cc}
{b_{2}-a_{2}(2+\frac{b_{2}r^{2}+b_{1}r+b_{0}}{a_{2}r^{2}+a_{1}r+a_{0}}),\,\,b_{1}+b_{2}r-(\frac{b_{2}r^{2}+b_{1}r+b_{0}}{a_{2}r^{2}+a_{1}r+a_{0}})(a_{1}+a_{2}r)-a_{1}}\\[2mm]
{a_{2}\,,\,a_{1}\,,\,\,a_{0}}
\end{array} & 
{x}
\end{array} \right).
\end{equation}

\begin{corollary}\label{coro3.1}
The eigenfunctions of the equation
\begin{multline}\label{eq:MainDif}
(x-r)\big(a_{2}x^{2} + a_{1}x + a_{0} \big) y_{n}''(x) + \big( b_{2} x^{2} + b_{1}x + b_{0} \big) y_{n}'(x) \\
-\Big(n\big(b_{2}+(n-1)a_{2}\big)(x-r)+c_{0}^{*}  \Big)y_{n}(x)=0\qquad \quad n\geq 1,
\end{multline}
where $ (-1)^{\frac{b_{2}r^{2}+b_{1}r+b_{0}}{a_{2}r^{2}+a_{1}r+a_{0}}}=1 $, are exceptional X$_{1}$-polynomials orthogonal with respect to the weight function \eqref{eq:10}. 
\end{corollary}
\medskip
\noindent
Let us make a contract here that the polynomial solution of equation \eqref{eq:MainDif} is indicated as
\begin{equation}\label{eq:denotingQ}
y_{n}(x)=Q_{n,r}\left(
\begin{array}{r|l}
\begin{array}{c}
{b_{2},b_{1},b_{0}}\\[2mm]
{a_{2},a_{1},a_{0}}\end{array} & 
{x}
\end{array} \right).
\end{equation}

By referring to the Pearson distributions family \eqref{eq:3}, we can now follow an inverse process and suppose that a simplified case of the weight function \eqref{eq:10} is given as 
\begin{equation}\label{eq:weighF}
\rho(x)=(x-r)^{\theta}\,W\left(
\begin{array}{r|l}
\begin{array}{cc}
{d^{*},\,\,e^{*}}\\[2mm]
{a\,,\,b\,,\,c}
\end{array} & 
{x}
\end{array} \right),
\end{equation} 
in which $ (-1)^{\theta}=1 $.
Then, by noting the equation \eqref{eq:MainDif} the unknown polynomials $ p_{2}(x) $ and $ q_{2}(x) $ of degree 2 in the differential equation 
\begin{equation}\label{eq:dif1}
(x-r)p_{2}(x)y''_{n}(x)+q_{2}(x)y'_{n}(x)-\big(\lambda_{n}(x-r)+c_{0}^{*}\big)y_{n}=0,
\end{equation}
can be directly derived by computing the logarithmic derivative of the function
\[ \frac{\rho(x)}{x-r}=(x-r)^{\theta-1}W\left(
\begin{array}{r|l}
\begin{array}{cc}
{d^{*},\,\,e^{*}}\\[2mm]
{a\,,\,b\,,\,c}
\end{array} & 
{x}
\end{array} \right)=(x-r)^{\theta-1}W(x), \]
as
\begin{align*}
\dfrac{\big((x-r)^{\theta-1}W(x)\big)'}{(x-r)^{\theta -1}W(x)}&=\dfrac{\theta -1}{x-r}+\dfrac{W'(x)}{W(x)}=\dfrac{\theta -1}{x-r}+\dfrac{d^{*}x+e^{*}}{ax^{2}+bx+c}\\[3mm]
&=\dfrac{\big(d^{*}+(\theta -1) a\big)x^{2}+\big(e^{*}-rd^{*}+(\theta -1) b\big)x-re^{*}+(\theta -1) c}{(x-r)(ax^{2}+bx+c)},
\end{align*}
and then equating the result with  
\[ \frac{q_{2}(x)-\big((x-r)p_{2}(x)\big)'}{(x-r)p_{2}(x)} ,\] 
so that we finally obtain
\begin{equation}\label{eq:2.13.1}
p_{2}(x)=ax^{2}+bx+c,
\end{equation}
and
\begin{equation}\label{eq:2.13.2}
q_{2}(x)=\big(d^{*}+(\theta +2)a\big)x^{2}+\big(e^{*}-r(d^{*}+2a)+(\theta+1)b\big)x+{\theta} c-r(e^{*}+b),
\end{equation}
provided that the roots of $ q_{2} $ are real.

Relations \eqref{eq:2.13.1} and \eqref{eq:2.13.2} show that the polynomial solution of equation \eqref{eq:dif1} can be represented in terms of the symbol \eqref{eq:denotingQ} as
\[y_{n}(x)=Q_{n,r}\left(
\begin{array}{r|l}
\begin{array}{c}
{d^{*}+(\theta +2)a,\,e^{*}-r(d^{*}+2a)+(\theta+1)b,\,{\theta} c-r(e^{*}+b)}\\[2mm]
{a,\,\,b,\,\,c}\end{array} & 
{x}
\end{array} \right),\] 
with the eigenvalue 
\[\lambda_{n}=n\big((n+1+\theta)a+d^{*}\big)\qquad n\geq 1.\]
Also, according to the Corollary \ref{coro1}, $ c_{0}^{*} $ in \eqref{eq:dif1} directly depends on the roots of $ q_{2}(x) $ in \eqref{eq:2.13.2} and is therefore computed as 
\begin{multline*}
c_{0}^{*}=2\theta (ar^{2}+br+c)\big(d^{*}+(\theta +2)a\big)\times\\
 \left( e^{*}+rd^{*}+(\theta +1)(2ra+b)
\mp {\Big(\big(e^{*}-r(d^{*}+2a)+(\theta+1)b\big)^{2}-4\big(d^{*}+(\theta +2)a\big)\big({\theta} c-r(e^{*}+b)\big)}\Big)^{\frac{1}{2}} \right) ^{-1}.
\end{multline*}

\medskip

As we observed, $\rho(x)$ was indeed the product of $(x-r)^{\theta} $ for 
\begin{equation}\label{eq:Theta}
\theta=\frac{b_{2}r^{2}+b_{1}r+b_{0}}{a_{2}r^{2}+a_{1}r+a_{0}},
\end{equation} 
and a special case of Pearson distributions family. This means that we can classify the exceptional X$ _{1} $-orthogonal polynomials into six main sequences.

\begin{corollary}\label{cor3.2}
By referring to table \ref{table1} and relation \eqref{eq:weighF}, there are totally six sequences of X$_{1}$-polynomials as follows:
\begin{enumerate}
\item Infinite X$_{1}$-Jacobi polynomials orthogonal with respect to the weight function
\begin{equation*}
\rho_{1}(x)=(x-r)^{\theta}(1-x)^{\alpha}(1+x)^{\beta}, \qquad (-1\leq x \leq 1).
\end{equation*}

\item Infinite X$_{1}$-Laguerre polynomials orthogonal with respect to the weight function
\[
\rho_{2}(x)=(x-r)^{\theta} x^{\alpha} \exp\left(-x \right), \qquad (0 \leq x < \infty).
\]

\item Infinite X$_{1}$-Hermite polynomials orthogonal with respect to the weight function
\[
\rho_{3}(x)=(x-r)^{\theta} \exp\left(-x^{2} \right), \qquad(-\infty<x<\infty).
\]

\item Finite X$_{1}$-polynomials orthogonal with respect to the weight function
\[
\rho_{4}(x)=(x-r)^{\theta} x^{q} (x+1)^{-(p+q)}, \qquad (0 \leq x < \infty).
\]

\item Finite X$_{1}$-polynomials orthogonal with respect to the weight function
\[
\rho_{5}(x)=(x-r)^{\theta} x^{-p} \exp \left( -\frac{1}{x} \right), \qquad (0\leq x < \infty).
\]

\item Finite X$_{1}$-polynomials orthogonal with respect to the weight function
\[
\rho_{6}(x)=(x-r)^{\theta} \left(1+x^{2} \right)^{-p} \exp(q \arctan x), \qquad (-\infty<x<\infty).
\]

\end{enumerate}
In all six above mentioned cases $r \in {\mathbb{R}}$ and $ \theta $ is a real parameter such that $ (-1)^{\theta}=1 $.
\end{corollary}

Note that for $\theta=-2$, the weight functions ( studied in \cite{MR2542180,MR2610341} and represented in \eqref{eq:Wxjacobi} and \eqref{eq:Wxlaguerre} ) are retrived for exceptional X$ _{1} $-Jacobi and X$ _{1} $-Laguerre polynomials when $ r=\frac{\beta+\alpha}{\beta-\alpha} $ and $ r=-\alpha $, respectively.

\medskip

\section{On the Series Solutions of Equation \eqref{eq:6} }
Let us reconsider equation \eqref{eq:6} in the form
\begin{multline}\label{eq.N1}
y_{n}''(x)+ \dfrac{ b_{2} x^{2} + b_{1}x + b_{0}  }{(x-r)\big(a_{2}x^{2} + a_{1}x+ a_{0} \big) }y_{n}'(x) 
-\dfrac{\lambda_{n}(x-r)+c_{0}^{*}  }{(x-r)\big(a_{2}x^{2} + a_{1}x +a_{0} \big) }y_{n}(x)=0,\\
\qquad  n \geq 1.
\end{multline}
Since
\begin{equation*}
\lim\limits_{x\rightarrow r}(x-r)\dfrac{ b_{2} x^{2} + b_{1}x + b_{0}  }{(x-r)\big(a_{2}x^{2} + a_{1}x + a_{0} \big) }=\dfrac{b_{2}r^{2}+b_{1}r+b_{0}}{a_{2}r^{2}+a_{1}r+a_{0}},
\end{equation*}
and
\begin{equation*}
\lim\limits_{x\rightarrow r}(x-r)^{2}\dfrac{\lambda_{n}(x-r)+c_{0}^{*}  }{(x-r)\big(a_{2}x^{2} +a_{1}x+ a_{0} \big) }=0,
\end{equation*}
the indicial equation corresponding to \eqref{eq.N1} is as
\begin{equation*}
t^{2}+\Big(\dfrac{b_{2}r^{2}+b_{1}r+b_{0}}{a_{2}r^{2}+a_{1}r+a_{0}}-1\Big)t=0.
\end{equation*}
By using the Frobenius method, one can obtain series solutions of equation \eqref{eq:6} when 
\[t_{1}=1-\dfrac{b_{2}r^{2}+b_{1}r+b_{0}}{a_{2}r^{2}+a_{1}r+a_{0}}=1-\theta,\] 
for different values of $ \theta $. 

If $ \theta\notin\mathbb{Z} $, the two basic solutions of \eqref{eq.N1} are respectively in the forms
\begin{equation*}
y_{n,1}(x)=\sum\limits_{k=0}^{\infty}C_{k}(x-r)^{k},\,\,\,\,\,C_{0}\neq 0,
\end{equation*}
and
\begin{equation*}
y_{n,2}(x)=(x-r)^{1-\theta}\sum\limits_{k=0}^{\infty}d_{k}(x-r)^{k},\,\,\,\,\,d_{0}\neq 0.
\end{equation*}
If $ \theta\in\mathbb{Z} $, three cases can occur for the basis solutions as follows: \\
\begin{equation*}
\theta=1\qquad\Rightarrow\quad \begin{cases}
y_{n,1}(x)=\sum\limits_{k=0}^{\infty}C_{k}(x-r)^{k},\,\,\,\,\,C_{0}\neq 0,\\
y_{n,2}(x)=y_{n,1}(x)\ln |x-r|+\sum\limits_{k=1}^{\infty}d_{k}(x-r)^{k},
\end{cases}
\end{equation*}
and
\begin{equation*}
\theta<1\qquad\Rightarrow\quad \begin{cases}
y_{n,1}(x)=(x-r)^{1-\theta}\sum\limits_{k=0}^{\infty}C_{k}(x-r)^{k},\,\,\,\,\,C_{0}\neq 0,\\
y_{n,2}(x)=wy_{n,1}(x)\ln |x-r|+\sum\limits_{k=0}^{\infty}d_{k}(x-r)^{k},\,\,\,\,\,d_{0}\neq 0,\quad w\in\mathbb{R},
\end{cases}
\end{equation*}
and finally
\begin{equation*}
\theta>1\qquad\Rightarrow\quad \begin{cases}
y_{n,1}(x)=\sum\limits_{k=0}^{\infty}C_{k}(x-r)^{k},\,\,\,\,\,C_{0}\neq 0,\\
y_{n,2}(x)=wy_{n,1}(x)\ln |x-r|+|x-r|^{1-\theta}\sum\limits_{k=0}^{\infty}d_{k}(x-r)^{k},\,\,\,\,\,d_{0}\neq 0,\quad w\in\mathbb{R}.
\end{cases}
\end{equation*}
Let us assume that 
\begin{equation}\label{eq.y.c1}
y_{n}(x)=\sum\limits_{k=0}^{\infty}C_{k}(x-r)^{k-\theta+1},
\end{equation}
for $ \theta\in\mathbb{Z}$ and $\theta<1 $. Since
\begin{equation*}
y'_{n}(x)=\sum\limits_{k=0}^{\infty}(k-\theta+1)C_{k}(x-r)^{k-\theta},\end{equation*}
and
\begin{equation*}
y''_{n}(x)=\sum\limits_{k=0}^{\infty}(k-\theta+1)(k-\theta)C_{k}(x-r)^{k-\theta-1},
\end{equation*}
substituting the above results in equation \eqref{eq:6} eventually leads to
\begin{multline}\label{eq.R1}
C_{k-1}\Big(a_{2}(k-\theta)(k-\theta-1)+b_{2}(k-\theta)-\lambda_{n}\Big)\\[3mm]
+C_{k}\Big((2a_{2}r+a_{1})(k-\theta+1)(k-\theta)+(2b_{2}r+b_{1})(k-\theta+1)-c_{0}^{*}\Big)\\[3mm]
+C_{k+1}\Big((a_{2}r^{2}+a_{1}r+a_{0})(k-\theta+2)(k-\theta+1)+(b_{2}r^{2}+b_{1}r+b_{0})(k-\theta+2)\Big)=0.
\end{multline}

In a similar way, for $ \theta\in\mathbb{Z}$ and $\theta\geq 1 $, or $ \theta\notin \mathbb{Z} $ the assumption
\begin{equation*}
y_{n}(x)=\sum\limits_{k=0}^{\infty}C_{k}(x-r)^{k},
\end{equation*}
leads to the same as recurrence relation \eqref{eq.R1} for $ \theta=1 $.
%\begin{multline}\label{eq.R.2}
%C_{k-1}\Big(a_{2}(k-1)(k-2)+b_{2}(k-1)-\lambda_{n}\Big)\\[3mm]
%+C_{k}\Big((2a_{2}r+a_{1})k(k-1)+(2b_{2}r+b_{1})k-c_{0}^{*}\Big)\\[3mm]
%+C_{k+1}\Big((a_{2}r^{2}+a_{1}r+a_{0})(k+1)k+(b_{2}r^{2}+b_{1}r+b_{0})(k+1)\Big)=0.
%\end{multline}

%In what follows, we focus on determining the polynomial solutions of equation \eqref{eq:MainDif} based on different values of the coefficients $ b_{0}, b_{1}, b_{2} $ and $ c_{0}^{*} $. 

\subsection{Some polynomials solutions of Equation \eqref{eq:6}}

According to corollary \ref{coro1}, since the coefficients of the polynomial $ B(x)=b_{2}x^{2}+b_{1}x+b_{0} $ in the main differential equation of \eqref{eq:6} or \eqref{eq:MainDif} have a significant role in determining the parameter $ c_{0}^{*} $ in relations \eqref{eq:2.2.1},  in this section we investigate six special cases of $ B(x) $ based on its roots and the real value $ r $ leading to particular cases of equation \eqref{eq:MainDif}.

\medskip
%\subsection*{Case 1}

First, suppose that $ b_{2}\neq 0 $ and $ r $ is a root of $ B(x) $. So,
\[b_{2}r^{2}+b_{1}r+b_{0}=0,\]
and relations \eqref{eq:2.2.1} are reduced to
\begin{equation}\label{case1}
\begin{cases}
2b_{2}r+b_{1}-c_{0}^{*}+\nu b_{2}=0,\\[2mm]
\nu c_{0}^{*}=0.
\end{cases}
\end{equation}
The relation $ \nu c_{0}^{*}=0 $ in \eqref{case1} gives three different cases as follows:\\
\\
$ \bullet $ {\bf Case 1.} $ \nu=0 $ and $ c_{0}^{*}=2b_{2}r+b_{1}=B'(r)\neq 0, $\\
\\
$ \bullet $ {\bf Case 2.} $ c_{0}^{*}=0 $ and $ \nu=-\dfrac{2b_{2}r+b_{1}}{b_{2}}=-\dfrac{B'(r)}{b_{2}}\neq 0, $\\
\\
$ \bullet $ {\bf Case 3.} $ c_{0}^{*}=0 $ and $ \nu=0, $ leading to $ B'(r)=0 $ which means that $ r $ is a multiple root of $ B(x) $.
% is of course in contradiction with the assumption.\\

\medskip

Second, suppose that $ b_{2}= 0 $ and $ b_{1}\neq 0 $. So, relations \eqref{eq:2.2.1} are reduced to
\begin{equation}\label{case2}
\begin{cases}
c_{0}^{*}=b_{1},\\[2mm]
\nu c_{0}^{*}=-(b_{1}r+b_{0}).
\end{cases}
\end{equation}
Now, if $ r $ is a root of $ B(x) $, we have $ \nu c_{0}^{*}=0 $ leading to \\
\\
$ \bullet $ {\bf Case 4.} $ c_{0}^{*}=b_{1}\neq 0 $ and $ \nu=0 $.\\
\\  
otherwise we get \\
\\
$ \bullet $ {\bf Case 5.} $ c_{0}^{*}=b_{1}\neq 0 $ and $ \nu=-\dfrac{b_{1}r+b_{0}}{c_{0}^{*}} =-\dfrac{B(r)}{b_{1}}\neq 0$.

%\medskip

Finally, suppose that $ b_{2}=b_{1}=0 $ so that relations \eqref{eq:2.2.1} are reduced to
\begin{equation}\label{case3}
\begin{cases}
c_{0}^{*}=0,\\[2mm]
b_{0}+\nu c_{0}^{*}=0, 
\end{cases}
\end{equation}
which yield $ b_{0}=0 $ leading to $ B(x)\equiv 0 $.  Therefore we get\\
\\
$ \bullet $ {\bf Case 6.} $ c_{0}^{*}= 0 $ and $ \nu $ is arbitrary.

\medskip

%
%
%Note that according to \eqref{eq:Theta}, in all above cases other than case 5, we have $ \theta=0 $. So, it is expected that under those circumstances the polynomial solutions of \eqref{eq:MainDif} will be the same as classical orthogonal polynomials. 
%
%
\medskip

\medskip

Under the conditions stated in Case 1, the differential equation \eqref{eq:MainDif} reads as
\begin{equation}\label{eq:4.4.1}
\big(a_{2}x^{2} + a_{1}x + a_{0} \big) y_{n}''(x) + \big(b_{2}x +  b_{2}r+b_{1}  \big) y_{n}'(x) 
-\Big(n\big(b_{2}+(n-1)a_{2}\big)+\dfrac{2b_{2}r+b_{1}}{x-r}  \Big)y_{n}(x)=0,
\end{equation}
for $  n\geq 1 $, whose solutions belong to the polynomial space
\[\Pi_{n,r,0}=\text{span}\Big\lbrace (x-r),(x-r)^{2},\ldots, (x-r)^{n}\Big\rbrace.\]
Without loss of generality, for $ r=0 $ the above equation takes the form
\begin{equation*}
(a_{2}x^{2} + a_{1}x + a_{0} ) y_{n}''(x) + ( b_{2}x + b_{1}  ) y_{n}'(x) 
-\Big(n\big(b_{2}+(n-1)a_{2}\big)+\dfrac{b_{1}}{x}  \Big)y_{n}(x)=0,
\end{equation*}
and for $ a_{0}=0 $ reads as
\begin{equation}\label{eq:special case dif}
x^{2}(a_{2}x+a_{1})y_{n}''(x) +x( b_{2}x + b_{1}  ) y_{n}'(x) 
-\big(n\big(b_{2}+(n-1)a_{2}\big) x+b_{1}  \big)y_{n}(x)=0.
\end{equation}
Note in this case that we have $ \theta=0 $ and for $ \lambda_{n}=n\big(b_{2}+(n-1)a_{2}\big) $, the coefficients $ C_{k} $ in \eqref{eq.R1} are recursively given by
\[C_{k}=\dfrac{\lambda_{n}-\lambda_{k}}{k\big(b_{1}+(k+1)a_{1}\big)}C_{k-1}\qquad \quad k=1,2,\ldots.\]
Setting $ C_{0}=1 $, it can be easily observed that 
\[C_{k}=\dfrac{1}{k!}\,\dfrac{\prod\limits_{j=0}^{k-1}\lambda_{n}-\lambda_{j+1}}{\prod\limits_{j=0}^{k-1}b_{1}+(j+2)a_{1}}\qquad\text{for}\qquad k=0,1,\ldots, n-1,\]
and $ C_{k}=0 $ for $ k\geq n $. So, by assuming
\[p_{n-1}(x)=\sum\limits_{k=0}^{n-1}C_{k}x^{k},\] 
and noting \eqref{eq.y.c1} we obtain $ y_{n}(x)=x\,p_{n-1}(x) $ as a polynomial solution of equation \eqref{eq:special case dif}. In this sense, replacing $ y'=p_{n-1}+xp'_{n-1} $ and $ y''=2p'_{n-1}+xp''_{n-1} $ in \eqref{eq:special case dif} yields the differential equation
\[(a_{2}x^{2}+a_{1}x)p''_{n-1}+\big((2a_{2}+b_{2})x+2a_{1}+b_{1}\big)p'_{n-1}-(n-1)(na_{2}+b_{1})p_{n-1}=0,\]
which has a polynomial solution of type \eqref{Dissertation.EQ.6.8} as
\[p_{n-1}(x)={P_{n-1}}\left( {\left. {\begin{array}{*{20}{c}}
{\begin{array}{*{20}{c}}
b_{2}+2a_{2},&b_{1}+2a_{1}
\end{array}}\\
{a_{2},\,\,\,\,\,\,a_{1},\,\,\,\,\,\,0}
\end{array}} \right|x} \right),\]
orthogonal with respect to a weight function of type \eqref{eq:3} as
\[W\left( {\left. {\begin{array}{*{20}{c}}
{\begin{array}{*{20}{c}}
b_{2},&a_{1}+b_{1}
\end{array}}\\
{a_{2},\,\,\,\,\,\,a_{1},\,\,\,\,\,\,0}
\end{array}} \right|x} \right)=x^{1+\frac{b_{1}}{a_{2}}}\big(a_{2}x+a_{1}\big)^{\frac{b_{2}}{a_{2}}-\frac{b_{1}}{a_{1}}-1}.\]

\begin{corollary}
The polynomial solution of the differential equation \eqref{eq:4.4.1}
is given by 
\[y_{n}(x)=(x-r){P_{n-1}}\left( {\left. {\begin{array}{*{20}{c}}
{\begin{array}{*{20}{c}}
b_{2}+2a_{2},&b_{1}+2a_{1}
\end{array}}\\
{a_{2},\,\,\,\,\,\,a_{1},\,\,\,\,\,\,a_{0}}
\end{array}} \right|x-r} \right).\]
This means that
\[{\bar{Q}_{n,r}}\left( {\left. {\begin{array}{*{20}{c}}
{\begin{array}{*{20}{c}}
b_{2},&b_{1},&-r(b_{2}r+b_{1})
\end{array}}\\
{a_{2},\, a_{1},\, a_{0}}
\end{array}} \right|x} \right)=
(x-r){\bar{P}_{n-1}}\left( {\left. {\begin{array}{*{20}{c}}
{\begin{array}{*{20}{c}}
b_{2}+2a_{2},&b_{1}+2a_{1}
\end{array}}\\
{a_{2},\,\,\,\,\,\,a_{1},\,\,\,\,\,\,a_{0}}
\end{array}} \right|x-r} \right).\]
\end{corollary}

Note that for $ b_{2}=0 $ in the above corollary the Case 4 is retrived.

\medskip

For cases 2, 3 and 6, where $ \theta $ is equal to zero, the differential equation \eqref{eq:MainDif} respectively reads as
\begin{equation}\label{eq1n}
\big(a_{2}x^{2} + a_{1}x + a_{0} \big) y_{n}''(x) + \big(b_{2}x +  b_{2}r+b_{1}  \big) y_{n}'(x) 
-n\big(b_{2}+(n-1)a_{2}\big) y_{n}(x)=0\qquad n\geq 1,
\end{equation}

\begin{equation}\label{eq2n}
\big(a_{2}x^{2} + a_{1}x + a_{0} \big) y_{n}''(x) + b_{2}(x-r) y_{n}'(x) 
-n\big(b_{2}+(n-1)a_{2}\big) y_{n}(x)=0\qquad n\geq 1,
\end{equation}
and
\begin{equation}\label{eq3n}
\big(a_{2}x^{2} + a_{1}x + a_{0} \big) y_{n}''(x) 
-n(n-1)a_{2}\, y_{n}(x)=0\qquad n\geq 1,
\end{equation}
which are all special cases of equation \eqref{EQ.hypDiff} and their polynomial solutions belong to the spaces
\[\text{span}\Big\lbrace \big(x+r+\dfrac{b_{1}}{b_{2}}\big),(x-r)^{2},\ldots, (x-r)^{n}\Big\rbrace,\]
$ \Pi_{n,r,0} $ and $ \Pi_{n,r,\nu} $ where $ \nu $ is an arbitrary value. 
\begin{corollary}
The monic polynomial solutions of equations \eqref{eq1n}-\eqref{eq3n} can be respectively denoted by
\begin{equation*}\label{eq:S1C2}
{\bar{Q}_{n,r}}\left( {\left. {\begin{array}{*{20}{c}}
{\begin{array}{*{20}{c}}
b_{2},&b_{1},&-r(rb_{2}+b_{1})
\end{array}}\\
{a_{2},\, a_{1},\, a_{0}}
\end{array}} \right|x} \right)={\bar{P}_n}\left( {\left. {\begin{array}{*{20}{c}}
{\begin{array}{*{20}{c}}
b_{2},& b_{2}r+b_{1}
\end{array}}\\
{a_{2},\,\,\,\,\,\,a_{1},\,\,\,\,\,\,a_{0}}
\end{array}} \right|x-r} \right),
\end{equation*}
\begin{equation*}
{\bar{Q}_{n,r}}\left( {\left. {\begin{array}{*{20}{c}}
{\begin{array}{*{20}{c}}
b_{2},&-2rb_{2},&r^{2}b_{2}
\end{array}}\\
{a_{2},\, a_{1},\, a_{0}}
\end{array}} \right|x} \right)={\bar{P}_n}\left( {\left. {\begin{array}{*{20}{c}}
{\begin{array}{*{20}{c}}
b_{2},& -b_{2}r
\end{array}}\\
{a_{2},\,\,\,\,\,\,a_{1},\,\,\,\,\,\,a_{0}}
\end{array}} \right|x-r} \right),
\end{equation*}
and
\begin{equation*}
{\bar{Q}_{n,r}}\left( {\left. {\begin{array}{*{20}{c}}
{\begin{array}{*{20}{c}}
0,&0,&0
\end{array}}\\
{a_{2},\, a_{1},\, a_{0}}
\end{array}} \right|x} \right)={\bar{P}_n}\left( {\left. {\begin{array}{*{20}{c}}
{\begin{array}{*{20}{c}}
0,& 0
\end{array}}\\
{a_{2},\,\,a_{1},\,\,a_{0}}
\end{array}} \right|x-r} \right).
\end{equation*}
\end{corollary}

\medskip

Finally, let us consider Case 5 where the differential equation \eqref{eq:MainDif} is reduced to
\begin{equation}\label{eq:C5}
\big(a_{2}x^{2} + a_{1}x + a_{0} \big) y_{n}''(x) + \big(b_{1} +\dfrac{b_{1}r+b_{0}}{x-r}\big)  y_{n}'(x) 
-\Big(n(n-1)a_{2}+\dfrac{b_{1}}{x-r}  \Big)y_{n}(x)=0,
\end{equation}
with the polynomial space
\[\Pi_{n,r,\nu}=\text{span}\Big\lbrace \big(x+\dfrac{b_{0}}{b_{1}}\big),(x-r)^{2},\ldots, (x-r)^{n}\Big\rbrace.\]

Again, without loss of generality, for $ r=0 $ the differential equation \eqref{eq:C5} takes the form
\begin{equation}\label{eq:special case5 dif}
\big(a_{2}x^{2} + a_{1}x + a_{0} \big) y_{n}''(x) + \big(b_{1} +\dfrac{b_{0}}{x}\big)  y_{n}'(x) 
-\Big(n(n-1)a_{2}+\dfrac{b_{1}}{x}  \Big)y_{n}(x)=0.
\end{equation}
In this case, for $ b_{2}=0, c_{0}^{*}=b_{1} $ and $ r=0 $, the  relation \eqref{eq.R1} can be written as
\begin{multline}\label{eq.C5}
C_{k-1}a_{2}\big((k-\sigma)(k-\sigma-1)-n(n-1)\big)\\[3mm]
+C_{k}(k-\sigma)\big(a_{1}(k-\sigma+1)+b_{1}\big)\\[3mm]
+C_{k+1}(k-\sigma+2)\big(a_{0}(k-\sigma+1)+b_{0}\big)=0,
\end{multline}
where
\[\sigma=\begin{cases}
\theta\qquad \theta<1 ,\,\theta\in\mathbb{Z}\\[1mm]
1\qquad \theta\geq 1  ,\,\theta\in\mathbb{Z}\quad \text{or}\quad \theta\notin\mathbb{Z}.
\end{cases}\]
Clearly, it is not possible generally to obtain $ C_{k} $ explicitly. However, in particular cases, if $ a_{2}=0 $ then $ y(x)=1+\dfrac{b_{1}}{b_{0}}x $ and if $ a_{1}=b_{1}=0 $, it is in contradiction with hypothesis in Case 5.
Moreover, $ a_{0}=b_{0}=0 $ lead to the special case 4 for $ \sigma=0 $. Let us consider \eqref{eq.C5} for $ a_{0}=b_{0}=0 $ when $ \sigma=1 $, which occurs if $\theta\in\mathbb{N}$ or $ \theta\notin\mathbb{Z}$, as
\begin{equation*}
C_{k-1}a_{2}\big((k-1)(k-2)-n(n-1)\big)
+C_{k}(k-1)\big(a_{1}k+b_{1}\big)=0.
\end{equation*}
So, the coefficients $ C_{k} $ in $ y_{n}(x)=\sum\limits_{k=1}^{\infty}C_{k}x^{k} $
 are recursively given by
\[C_{k}=\dfrac{n(n-1)-(k-1)(k-2)}{(k-1)\big(a_{1}k+b_{1}\big)}a_{2}C_{k-1}\qquad\text{for} \quad k=2,3,\ldots.\]
Setting $ C_{1}=1 $, it can be easily observed that 
\[C_{k}=\dfrac{a_{2}^{k-1}}{(k-1)!}\,\dfrac{\prod\limits_{j=0}^{k-2}n(n-1)-j(j+1)}{\prod\limits_{j=0}^{k-2}b_{1}+(j+2)a_{1}}\qquad\text{for}\qquad k=1,2,3,\ldots, n,\]
and $ C_{k}=0 $ for $ k> n $.

\medskip

\section{Six Classes of Exceptional X$ _{1} $-Orthogonal Polynomials}

By noting the corollary \ref{cor3.2}, in this section we consider six special cases of the main equation \eqref{eq:MainDif} and study their properties. 

\subsection{On the differential equation of exceptional X$ _{1} $-Jacobi Polynomials}

As a generalization of Jacobi differential equation for $ \theta=0 $, consider the equation
\begin{multline}\label{diffX1Jacobi}
(x-r)(1-x^{2})y''_{n}(x)+\Big(-(\alpha+\beta+\theta+2)x^{2}+\big(\beta-\alpha+r(\alpha+\beta+2)\big)x+\theta-r(\beta-\alpha)\Big)y'_{n}(x)\\
+\Big(n(n+\alpha+\beta+\theta+1)(x-r)-c_{0}^{(P)}\Big)y_{n}(x)=0\qquad\quad n\geq 1,
\end{multline}
where $ r,\theta, \alpha,\beta $ are real parameters such that $ \alpha,\beta>-1 $, $ (-1)^{\theta}=1 $ and
\begin{multline*}
c_{0}^{(P)}={2\theta(1-r^{2})(\alpha+\beta+\theta+2)}
\Big( (\alpha+\beta+2\theta+2)r+\alpha-\beta \\
 \pm \sqrt{\big((\alpha+\beta+2)r+\beta-\alpha\big)^{2}+4(\alpha+\beta+\theta+2)\big(\theta-r(\beta-\alpha)\big)}\,\Big) ^{-1}.
\end{multline*}
The polynomial solution of equation \eqref{diffX1Jacobi}, i.e.
\[y_{n}(x)=P_{n,r,\theta}^{(\alpha,\beta)}(x)={Q}_{n,r}\left(\begin{array}{r|l}\begin{array}{c}
-(\alpha+\beta+\theta+2),\,\,\beta-\alpha+r(\alpha+\beta+2),\,\,\theta-r(\beta-\alpha)\\[3mm]
-1,\,\,0,\,\,1\end{array} & {x}\end{array} \right),\]
is orthogonal with respect to the weight function
\begin{equation*}
\rho_{1}(x;\,\alpha,\beta,\theta)=(x-r)^{\theta}(1-x)^{\alpha}(1+x)^{\beta}=(x-r)^{\theta}\,W\left(\begin{array}{r|l}\begin{array}{c}-\alpha-\beta,\,\,\beta-\alpha\\-1,\,\,0,\,\,1\end{array} & {x}\end{array} \right), 
\end{equation*}
on $ [-1,1] $. Also, for $ \theta=0 $, $ r=-1 $ or $ r=1 $ in \eqref{diffX1Jacobi}, $ c_{0}^{(P)}=0 $ and the weight function $ \rho_{1}(x;\,\alpha,\beta,\theta) $ will be a special case of the beta distribution. Hence, the solution of equation \eqref{diffX1Jacobi} will be the same as classical Jacobi polynomials. In fact, in each of these circumstances equation \eqref{diffX1Jacobi} reads as
\begin{equation*}
(1-x^{2})y''_{n}(x)+\big(-(\alpha+\beta+2)x+\beta-\alpha\big)y'_{n}(x)
+n(n+\alpha+\beta+1)y_{n}(x)=0,
\end{equation*}
for $ \theta=0 $ and
\begin{equation*}
(1-x^{2})y''_{n}(x)+\big(-(\alpha+\beta+\theta+2)x+\beta+\theta-\alpha\big)y'_{n}(x)
+n(n+\alpha+\beta+\theta+1)y_{n}(x)=0,
\end{equation*}
for $ r=-1 $ and
\begin{equation*}
(1-x^{2})y''_{n}(x)+\big(-(\alpha+\beta+\theta+2)x+\beta-\theta-\alpha\big)y'_{n}(x)
+n(n+\alpha+\beta+\theta+1)y_{n}(x)=0,
\end{equation*}
for $ r=1 $ with the following Jacobi type polynomial solutions
\begin{align*}
& {P}_{n,r,0}^{(\alpha,\beta)}(x)={P}_{n}^{(\alpha,\beta)}(x),\\[2mm]
& {P}_{n,-1,\theta}^{(\alpha,\beta)}(x)={P}_{n}^{(\alpha,\beta+\theta)}(x),
\end{align*}
and
\[{P}_{n,1,\theta}^{(\alpha,\beta)}(x)={P}_{n}^{(\alpha+\theta,\beta)}(x).\]

\subsection{On the differential equation of exceptional X$ _{1} $-Laguerre Polynomials}

As a generalization of Laguerre differential equation for $ \theta=0 $, consider the equation
\begin{multline}\label{diffX1Laguerre}
x(x-r) y''_{n}(x)+\Big(-x^{2}+(\alpha+r+\theta+1)x-r(\alpha+1)\Big)y'_{n}(x)
+\big(n(x-r)-c_{0}^{(L)}\big)y_{n}(x)=0\\
\qquad\quad n\geq 1,
\end{multline}
where $ r,\theta, \alpha $ are real parameters such that $ \alpha>-1 $, $ (-1)^{\theta}=1 $ and
\begin{equation*}
c_{0}^{(L)}={2r\theta}\left({r-\alpha-\theta-1\pm \sqrt{(r+\theta)^{2}+(\alpha+1)(\alpha+1+2\theta-2r)}}\right)^{-1}.
\end{equation*}
The polynomial solution of equation \eqref{diffX1Laguerre}, i.e.
\[{L}_{n,r,\theta}^{(\alpha)}(x)={Q}_{n,r}\left(\begin{array}{r|l}\begin{array}{c}
-1,\,\,\alpha+r+\theta+1,\,\,-r(\alpha+1) \\[3mm]
0,\,\,1,\,\,0\end{array} & {x}\end{array} \right),\]
is orthogonal with respect to the weight function
\begin{equation*}
\rho_{2}(x;\,\alpha,\theta)=(x-r)^{\theta}x^{\alpha}e^{-x}=(x-r)^{\theta}\,W\left(\begin{array}{r|l}\begin{array}{c}
-1,\,\,\alpha\\
0,\,\,1,\,\,0
\end{array} & {x}\end{array} \right), 
\end{equation*}
on $ [0,\infty) $. Also, for $ \theta=0 $ or $ r=0 $ in \eqref{diffX1Laguerre}, $ c_{0}^{(L)}=0 $ and the weight function $ \rho_{2}(x;\,\alpha,\theta) $ will be a special case of Gamma distribution. Hence, the solution of equation \eqref{diffX1Laguerre} will be the same as classical Laguerre polynomials. In fact, in each of these circumstances equation \eqref{diffX1Laguerre} reads as
\begin{equation*}
x y''_{n}(x)+(-x+\alpha+1)y'_{n}(x)
+ny_{n}(x)=0,
\end{equation*}
for $ \theta=0 $ and
\begin{equation*}
x y''_{n}(x)+(-x+\alpha+\theta+1)y'_{n}(x)
+ny_{n}(x)=0,
\end{equation*}
for $ r=0 $ with the following Laguerre type polynomial solutions
\[{L}_{n,r,0}^{(\alpha)}(x)={L}_{n}^{(\alpha)}(x),\]
and
\[{L}_{n,0,\theta}^{(\alpha)}(x)={L}_{n}^{(\alpha+\theta)}(x).\]

\subsection{On the differential equations of exceptional X$ _{1} $-Hermite Polynomials}

As a generalization of Hermite differential equation for $ \theta=0 $, consider the equation
\begin{equation}\label{diffX1Hermite}
(x-r) y''_{n}(x)+(-2x^{2}+2rx+\theta)y'_{n}(x)
+\Big(2n(x-r)-\dfrac{2\theta}{r\pm\sqrt{r^{2}+2\theta}}\Big)y_{n}(x)=0\,\,\,\qquad n\geq 1,
\end{equation}
where $ r,\theta $ are real parameters such that $ (-1)^{\theta}=1 $.

The polynomial solution of equation \eqref{diffX1Hermite}, i.e. 
\[{H}_{n,r,\theta}(x)={Q}_{n,r}\left(\begin{array}{r|l}\begin{array}{c}
-2,\,\,2r,\,\,\theta\\[3mm]
0,\,\,0,\,\,1\end{array} & {x}\end{array} \right),\]
is orthogonal with respect to the weight function
\begin{equation*}
\rho_{3}(x; \,\theta)=(x-r)^{\theta}e^{-x^{2}}=(x-r)^{\theta}\,W\left(\begin{array}{r|l}\begin{array}{c}
-2,\,\,0\\
0,\,\,0,\,\,1
\end{array} & {x}\end{array} \right), 
\end{equation*}
on $ (-\infty,\infty) $. Also, for $ \theta=0 $, $ \rho_{3}(x; \,\theta) $ is the same as normal distribution and the solution of equation \eqref{diffX1Hermite} is the classical Hermite polynomials. In fact, in this case, equation \eqref{diffX1Hermite} reads as
\begin{equation*}
y''_{n}(x)-2xy'_{n}(x)+2ny_{n}(x)=0,
\end{equation*}
with the usual Hermite polynomial solution
\[{H}_{n,r,0}(x)={H}_{n}(x).\]

\subsection{First Finite Sequence of Exceptional X$ _{1} $-Orthogonal Polynomials}

Consider the following differential equation
\begin{multline}\label{X1F1Diff}
x(x-r)(x+1)y''_{n}(x)+\Big((\theta+2-p)x^{2}+\big(q+\theta+1+r(p-2)\big)x-r(q+1)\Big)y'_{n}(x)\\
-\Big(n(n+1+\theta-p)(x-r)+c_{0}^{(M)}\Big)y_{n}(x)=0\qquad\quad n\geq 1,
\end{multline}
where $ r,\theta $ are real parameters such that $ (-1)^{\theta}=1 $ and
\[c_{0}^{(M)}=\dfrac{2\theta r(r+1)(p-\theta-2)}{rp-q-(\theta+1)(2r+1)\pm\Big(\big(q+\theta+1+r(p-2)\big)^{2}+4r(q+1)(\theta+2-p)\Big)^{\frac{1}{2}}}.\]
The polynomial solution of equation \eqref{X1F1Diff}, i.e. 
\[{M}_{n,r,\theta}^{(p,q)}(x)={Q}_{n,r}\left(\begin{array}{r|l}\begin{array}{c}
\theta+2-p,\,\,q+\theta+1+r(p-2),\,\,-r(q+1)\\[3mm]
1,\,\,1,\,\,0\end{array} & {x}\end{array} \right),\]
is finitely orthogonal with respect to the weight function
\begin{equation*}
\rho_{4}(x; p,q,\theta)=(x-r)^{\theta}x^{q}(x+1)^{-(p+q)}=(x-r)^{\theta}\,W\left(\begin{array}{r|l}\begin{array}{c}
-p,\,\,q\\
1,\,\,1,\,\,0\end{array} & {x}\end{array} \right), 
\end{equation*}
on $ [0,\infty) $ if and only if
\[ p>2\lbrace\max{n}\rbrace+\theta+1 \qquad\text{and}\qquad  q>-1 .\]
In other words, if the self-adjoint form of equation \eqref{X1F1Diff} is written as
\begin{multline}\label{selfF1n}
\Big((x-r)^{\theta}x^{q+1}(x+1)^{1-(p+q)}y'_{n}(x)\Big)'\\=(x-r)^{\theta-1}x^{q}(x+1)^{-(p+q)}\Big(n(n+1+\theta-p)(x-r)+c_{0}^{(M)}\Big)y_{n}(x),
\end{multline}
and for the index $ m $ as
\begin{multline}\label{selfF1m}
\Big((x-r)^{\theta}x^{q+1}(x+1)^{1-(p+q)}y'_{m}(x)\Big)'\\=(x-r)^{\theta-1}x^{q}(x+1)^{-(p+q)}\Big(m(m+1+\theta-p)(x-r)+c_{0}^{(M)}\Big)y_{m}(x),
\end{multline}
then multiplying by $ y_{m}(x) $ and $ y_{n}(x) $ in relations \eqref{selfF1n} and \eqref{selfF1m} respectively and subtracting them and finally integrating from both sides gives
\begin{multline}\label{EQ.SLF1}
\Big[ (x-r)^{\theta}x^{q+1}(x+1)^{1-(p+q)}\big(y'_{n}(x)y_{m}(x)-y'_{m}(x)y_{n}(x)\big)\Big]_{0}^{\infty}\\
=\big(n(n+1+\theta-p)-m(m+1+\theta-p)\big)\int_{0}^{\infty}(x-r)^{\theta}x^{q}(x+1)^{-(p+q)}y_{n}(x)y_{m}(x)\,dx.
\end{multline}
Now, since 
\[\text{max deg}\,\,\{ {y'_n}(x){y_m}(x) - {y'_m}(x){y_n}(x)\}  = m + n - 1,\]
if 
\[ q>-1 \quad\text{and} \quad p>2N+\theta+1 \quad\text{for} \quad N=\max\lbrace m,n\rbrace,\] 
the left hand side of \eqref{EQ.SLF1} tends to zero and for $ m,n\geq 1 $ we get
\begin{multline*}
\int_0^\infty  {\,\,\frac{{(x-r)^{\theta}{x^q}}}{{{{(x + 1)}^{p + q}}}}\,M_{n,r,\theta}^{(p,q)}(x)\,M_{m,r,\theta}^{(p,q)}(x)\,dx}  = \,0\,\,\,\\
 \Leftrightarrow \,\,\,
m \ne n,\,\,N=\max\lbrace m,n\rbrace<\dfrac{p-1-\theta}{2},\,\,q >  - 1\,\,\,\,\text{and}\,\,(-1)^{\theta}=1.
\end{multline*}

Note that for $ \theta=0 $, $ r=-1 $ or $ r=0 $, $ \rho_{4}(x; p,q, \theta) $ would be a special case of  F-Fisher distribution. Indeed, in each of these circumstances $ c_{0}^{(M)}=0 $ and equation \eqref{X1F1Diff} reads as
\begin{equation*}
x(x+1)y''_{n}(x)+\big((2-p)x+q+1\big)y'_{n}(x)
-n(n+1-p)y_{n}(x)=0,
\end{equation*}
for $ \theta=0 $ and
\begin{equation*}
x(x+1)y''_{n}(x)+\big((\theta+2-p)x+q+1\big)y'_{n}(x)
-n(n+1+\theta-p)y_{n}(x)=0,
\end{equation*}
for $ r=-1 $ and
\begin{equation*}
x(x+1)y''_{n}(x)+\big((\theta+2-p)x+q+\theta+1\big)y'_{n}(x)
-n(n+1+\theta-p)y_{n}(x)=0,
\end{equation*}
for $ r=0 $ with the following polynomial solutions
\begin{align*}
& {M}_{n,r,0}^{(p,q)}(x)={M}_{n}^{(p,q)}(x),\\[2mm]
& {M}_{n,-1,\theta}^{(p,q)}(x)={M}_{n}^{(p-\theta,q)}(x),
\end{align*}
and
\[{M}_{n,0,\theta}^{(p,q)}(x)={M}_{n}^{(p-\theta,q+\theta)}(x).\]

\medskip

\subsection{Second Finite Sequence of Exceptional X$ _{1} $-Orthogonal Polynomials}

Consider the following differential equation
\begin{multline}\label{X1F2Diff}
(x-r)x^{2} y''_{n}(x)+\Big((\theta+2-p)x^{2}+\big(1+r(p-2)\big)x-r\Big)y'_{n}(x)\\
-\Big(n(n+1+\theta-p)(x-r)+c_{0}^{(N)}\Big)y_{n}(x)=0\qquad\quad n\geq 1,
\end{multline}
where $ r,\theta $ are real parameters such that $ (-1)^{\theta}=1 $ and
\[c_{0}^{(N)}=\dfrac{2\theta r^{2}(p-\theta-2)}{r\big(p-2(\theta+1)\big)-1\pm \Big(\big(1+r(p-2)\big)^{2}+4r(\theta+2-p)\Big)^{\frac{1}{2}}}.\]
The polynomial solution of equation \eqref{X1F2Diff}, i.e.
\[{N}_{n,r,\theta}^{(p)}(x)={Q}_{n,r}\left(\begin{array}{r|l}\begin{array}{c}
\theta+2-p,\,\,1+r(p-2),\,\,-r\\[3mm]
1,\,\,0,\,\,0\end{array} & {x}\end{array} \right),\]
is finitely orthogonal with respect to the weight function
\begin{equation*}
\rho_{5}(x; p,\theta)=(x-r)^{\theta}x^{-p}e^{-\frac{1}{x}}=(x-r)^{\theta}\,W\left(\begin{array}{r|l}\begin{array}{c}
-p,\,\,1\\
1,\,\,0,\,\,0\end{array} & {x}\end{array} \right), 
\end{equation*}
on $ [0,\infty) $ if and only if $ p>2\lbrace\max{n}\rbrace+\theta+1 $, because if the self-adjoint form of equation \eqref{X1F2Diff} is written as
\begin{multline}\label{selfF2n}
\Big((x-r)^{\theta}x^{-p+2}e^{-\frac{1}{x}}y'_{n}(x)\Big)'\\=(x-r)^{\theta-1}x^{-p}e^{-\frac{1}{x}}\Big(n(n+1+\theta-p)(x-r)+c_{0}^{(N)}\Big)y_{n}(x),
\end{multline}
and for the index $ m $ as
\begin{multline}\label{selfF2m}
\Big((x-r)^{\theta}x^{-p+2}e^{-\frac{1}{x}}y'_{m}(x)\Big)'\\=(x-r)^{\theta-1}x^{-p}e^{-\frac{1}{x}}\Big(m(m+1+\theta-p)(x-r)+c_{0}^{(N)}\Big)y_{m}(x),
\end{multline}
then multiplying by $ y_{m}(x) $ and $ y_{n}(x) $ in relations \eqref{selfF2n} and \eqref{selfF2m} respectively and subtracting them and finally integrating from both sides gives
\begin{multline}\label{EQ.SLF2}
\Big[ (x-r)^{\theta}x^{-p+2}e^{-\frac{1}{x}}\big(y'_{n}(x)y_{m}(x)-y'_{m}(x)y_{n}(x)\big)\Big]_{0}^{\infty}\\
=\big(n(n+1+\theta-p)-m(m+1+\theta-p)\big)\int_{0}^{\infty}(x-r)^{\theta}x^{-p}e^{-\frac{1}{x}}y_{n}(x)y_{m}(x)\,dx.
\end{multline}
Now, if 
\[ p>2N+\theta+1 \quad\text{for} \quad N=\max\lbrace m,n\rbrace,\] 
the left hand side of \eqref{EQ.SLF2} tends to zero and for $ m,n\geq 1 $ we get
\begin{multline*}
\int_{0}^{\infty}(x-r)^{\theta}{x^{-p}}e^{-\frac{1}{x}}\,N_{n,r,\theta}^{(p)}(x)\,N_{m,r,\theta}^{(p)}(x)\,dx = \,0\,\,\, \\
\Leftrightarrow \,\,\,
m \ne n,\,\,N=\max\lbrace m,n\rbrace < \dfrac{p-\theta-1}{2}\,\,\,\text{and}\,\,\,(-1)^{\theta}=1.
\end{multline*}

Note that for $ \theta=0 $ or $ r=0 $, $ \rho_{5}(x; p,\theta) $ would be a special case of  inverse Gamma distribution. Indeed, in each of these circumstances $ c_{0}^{(N)}=0 $ and equation \eqref{X1F2Diff} reads as
\begin{equation*}
x^{2} y''_{n}(x)+\big((2-p)x+1\big)y'_{n}(x)
-n(n+1-p)y_{n}(x)=0,
\end{equation*}
for $ \theta=0 $ and
\begin{equation*}
x^{2} y''_{n}(x)+\big((\theta+2-p)x+1\big)y'_{n}(x)
-n(n+1+\theta-p)y_{n}(x)=0,
\end{equation*}
for $ r=0 $ with the following polynomial solutions
\[{N}_{n,r,0}^{(p)}(x)={N}_{n}^{(p)}(x),\]
and
\[ {N}_{n,0,\theta}^{(p)}(x)=\bar{N}_{n}^{(p-\theta)}(x).\]

\subsection{Third Finite Sequence of Exceptional X$ _{1} $-Orthogonal Polynomials}

Consider the following differential equation
\begin{multline}\label{X1F3Diff}
(x-r)(1+x^{2}) y''_{n}(x)+\Big((\theta+2-2p)x^{2}+\big(q+2r(p-1)\big)x+\theta-rq\Big)y'_{n}(x)\\
-\Big(n(n+1+\theta-2p)(x-r)+c_{0}^{(J)}\Big)y_{n}(x)=0\qquad\quad n\geq 1,
\end{multline}
where $ r,\theta $ are real parameters such that $ (-1)^{\theta}=1 $ and
\[c_{0}^{(J)}=\dfrac{2\theta (r^{2}+1)(2p-\theta-2)}{2r(p-\theta-1)-q\pm\Big(\big(q+2r(p-1)\big)^{2}-4(\theta+2-2p)(\theta-rq)\Big)^{\frac{1}{2}}}.\]
The polynomial solution of equation \eqref{X1F3Diff}, i.e.
\[{J}_{n,r,\theta}^{(p,q)}(x)={Q}_{n,r}\left(\begin{array}{r|l}\begin{array}{c}
\theta+2-2p,\,\,q+2r(p-1),\,\,\theta -rq\\[3mm]
1,\,\,0,\,\,1\end{array} & {x}\end{array} \right),\]
is finitely orthogonal with respect to the weight function
\begin{equation*}
\rho_{6}(x; p,q,\theta)=(x-r)^{\theta}(1+x^{2})^{-p}\exp(q\arctan x)=(x-r)^{\theta}\,W\left(\begin{array}{r|l}\begin{array}{c}
-2p,\,\,q\\
1,\,\,0,\,\,1\end{array} & {x}\end{array} \right), 
\end{equation*}
on $ (-\infty,\infty) $ if and only if $ p>\lbrace\max{n}\rbrace+\frac{\theta+1}{2} $, because if the self-adjoint form of equation \eqref{X1F3Diff} is written as
\begin{multline}\label{selfF3n}
\Big((x-r)^{\theta}{\left( {{1+x^2} } \right)^{1 - p}}\exp (q\arctan x){y'_n}(x)\Big)^\prime 
\\
 = (x-r)^{\theta-1}{\left( 1+x^2 \right)^{ - p}}\exp (q\arctan x)\Big(n(n + 1+\theta - 2p)(x-r)+c_{0}^{(J)}\Big){y_n}(x),
\end{multline}
and for the index $ m $ as
\begin{multline}\label{selfF3m}
\Big((x-r)^{\theta}{\left( {{1+x^2} } \right)^{1 - p}}\exp (q\arctan x){y'_m}(x)\Big)^\prime 
\\
 = (x-r)^{\theta-1}{\left( 1+x^2 \right)^{ - p}}\exp (q\arctan x)\Big(m(m + 1+\theta - 2p)(x-r)+c_{0}^{(J)}\Big){y_m}(x),
\end{multline}
then multiplying by $ y_{m}(x) $ and $ y_{n}(x)  $ in relations \eqref{selfF3n} and \eqref{selfF3m} respectively and subtracting them and finally integrating from both sides gives
\begin{multline}\label{EQ.SLF3}
\left[(x-r)^{\theta}{\left( {{1+x^2} } \right)^{1 - p}}\exp (q\arctan x)\big({y'_n}(x){y_m}(x) - {y'_m}(x)\,{y_n}(x)\big)\right]_{ - \infty }^\infty  \\
= \big( n(n+1+\theta-2p)-m(m+1+\theta-2p) \big)\int_{ - \infty }^\infty  {(x-r)^{\theta}{{\left( 1+x^2 \right)}^{ - p}}\exp (q\arctan x)\,J_{n,r,\theta}^{(p,q)}(x)J_{m,r,\theta}^{(p,q)}(x)\,dx}.
\end{multline}
Now, if 
\[ p > N + \frac{\theta+1}{2}\qquad\text{for}\qquad N = \max \{ m,n\},\] 
the left hand side of \eqref{EQ.SLF3} tends to zero and for $ m,n\geq 1 $ we have 
\begin{multline*}
\int_{ - \infty }^\infty  (x-r)^{\theta}{(1+x^2)^{ - p}}\exp (q\arctan x)\,J_{n,r,\theta}^{(p,q)}(x)J_{m,r,\theta}^{(p,q)}(x)\,dx  = 0\\
 \Leftrightarrow \quad
m \ne n\,,\,N= \max \{ m,n\}<p-\frac{\theta+1}{2}\,\,\text{and}\,\,\,(-1)^{\theta}=1 .
\end{multline*}

For $ \theta=0 $, $ \rho_{6}(x;p,q,\theta) $ is reduced to the generalized T-Student distribution. In this case $ c_{0}^{(J)}=0 $ and equation \eqref{X1F3Diff} reads as
\begin{equation*}
(1+x^{2}) y''_{n}(x)+\big(2(1-p)x+q\big)y'_{n}(x)
-n(n+1-2p)y_{n}(x)=0,
\end{equation*}
with the polynomial solution
\[{J}_{n,r,0}^{(p,q)}(x)={J}_{n}^{(p,q)}(x).\]

\section*{Acknowledgments}
This work has been supported by the {\em{Alexander von Humboldt Foundation}} under the grant number: Ref 3.4 - IRN - 1128637 - GF-E.


\begin{thebibliography}{99}

\bibitem{AlSalam}
W. Al-Salam, Characterization theorems for orthogonal polynomial. Orthogonal Polynomials: Theory and Practice, \emph{NATO ASI Series}, Vol. 294. Kluwer Acad., Dordrecht, Netherlands, (1990) pp. 1–24. 

\bibitem{Andrews}
G. E. Andrews, R. Askey and R. Roy: {\em Special Functions. Encyclopedia of Mathematics and its Applications} 71, Cambridge University Press, Cambridge 1999.




\bibitem{Arfken}
G. Arfken, {\em Mathematical methods for physicists}, Academic Press Inc 1985.


\bibitem{MR1545034}
S.~Bochner.
\newblock \"{U}ber {S}turm-{L}iouvillesche {P}olynomsysteme.
\newblock {\em Math. Z.}, 29(1):730--736, 1929.

\bibitem{MR0481884}
T.~S. Chihara.
\newblock {\em An introduction to orthogonal polynomials}.
\newblock Gordon and Breach Science Publishers, New York, 1978.

\bibitem{MR2905628}
D.~Dutta and P.~Roy.
\newblock Darboux transformation, exceptional orthogonal polynomials and information theoretic measures of uncertainty.
\newblock In {\em Algebraic aspects of {D}arboux transformations, quantum integrable systems and supersymmetric quantum mechanics}, volume 563 of {\em Contemp. Math.}, pages 33--49. Amer. Math. Soc., Providence, RI, 2012.

\bibitem{Freund}
J. E. Freund, R. E. Walpole, \emph{Mathematical Statistics}, Prentice-Hall. Inc., New Jersey, USA. 1971.

\bibitem{Ferrero}
M. Garc{\'ia}-Ferrero, D. G{\'o}mez-Ullate, R. Milson, A Bochner type characterization theorem for exceptional orthogonal polynomials, \emph{J. Math. Anal. Appl.}, 472(1), (2019), 584-626.


\bibitem{MR2542180}
D. G{{\'o}}mez-Ullate, N. Kamran, R. Milson.
\newblock An extended class of orthogonal polynomials defined by a {S}turm-{L}iouville problem.
\newblock {\em J. Math. Anal. Appl.}, 359(1):352--367, 2009.

\bibitem{MR2610341}
D. G{{\'o}}mez-Ullate, N. Kamran, R. Milson,
\newblock An extension of {B}ochner's problem: exceptional invariant subspaces.
\newblock {\em J. Approx. Theory}, 162(5):987--1006, 2010.

\bibitem{MR2771725}
Choon-Lin Ho.
\newblock Dirac(-{P}auli), {F}okker-{P}lanck equations and exceptional {L}aguerre polynomials.
\newblock {\em Ann. Physics}, 326(4):797--807, 2011.

\bibitem{MR2191786}
Mourad E.~H. Ismail.
\newblock {\em Classical and quantum orthogonal polynomials in one variable}.
\newblock Cambridge University Press, Cambridge, 2005.

\bibitem{Kotz}
Norman~L. Johnson, Samuel Kotz, and N.~Balakrishnan.
\newblock {\em Continuous univariate distributions (2nd edition)}.
\newblock John Wiley \& Sons Inc., 1994.

\bibitem{MR2656096}
R.~Koekoek, P.~A. Lesky, and R.~F. Swarttouw.
\newblock {\em Hypergeometric orthogonal polynomials and their {$q$}-analogues}.
\newblock Springer Monographs in Mathematics. Springer-Verlag, Berlin, 2010.

\bibitem{MR2246501}
W. Koepf, M. Masjed-Jamei.
\newblock A generic polynomial solution for the differential equation of hypergeometric type and six sequences of orthogonal polynomials related to it.
\newblock {\em Integral Transforms Spec. Funct.}, 17(8):559--576, 2006.

\bibitem{MR2053407}
M. Masjed-Jamei.
\newblock Classical orthogonal polynomials with weight function {$((ax+b)^2+(cx+d)^2)^{-p}\exp(q\,{\rm Arctg}((ax+b)/(cx+d)))$},
  {$x\in(-\infty,\infty)$} and a generalization of {$T$} and {$F$} distributions.
\newblock {\em Integral Transforms Spec. Funct.}, 15(2):137--153, 2004.

\bibitem{MR1915513}
M. Masjed-Jamei.
\newblock Three finite classes of hypergeometric orthogonal polynomials and their application in functions approximation.
\newblock {\em Integral Transforms Spec. Funct.}, 13(2):169--191, 2002.

\bibitem{theBook}
M. Masjed-Jamei, {\em Special functions and generalized Sturm-Liouville problems}, Springer Nature, Birkh\"{a}user, 2020.

\bibitem{MMS}
M. Masjed-Jamei, Z. Moalemi, N. Saad. 
\newblock Incomplete symmetric orthogonal polynomials of finite type generated by a generalized Sturm-Liouville theorem.
\newblock {\em J. Math. Phys.}, 61(2):023501, 2020.

\bibitem{MR1149380}
A.~F. Nikiforov, S.~K. Suslov, and V.~B. Uvarov.
\newblock {\em Classical orthogonal polynomials of a discrete variable}.
\newblock Springer Series in Computational Physics. Springer-Verlag, Berlin, 1991.

\bibitem{Niki}
A. F. Nikiforov, V. B. Uvarov, \emph{Special Functions of Mathematical Physics}, Basel-Boston: Birkh\"{a}user, 1988. 


\bibitem{MR2569488}
Satoru Odake and Ryu Sasaki.
\newblock Infinitely many shape invariant potentials and new orthogonal polynomials.
\newblock {\em Phys. Lett. B}, 679(4):414--417, 2009.

\bibitem{MR2588057}
Satoru Odake and Ryu Sasaki.
\newblock Another set of infinitely many exceptional {$(X_\ell)$} {L}aguerre polynomials.
\newblock {\em Phys. Lett. B}, 684(2-3):173--176, 2010.

\bibitem{MR2439200}
C.~Quesne.
\newblock Exceptional orthogonal polynomials, exactly solvable potentials and supersymmetry.
\newblock {\em J. Phys. A}, 41(39):392001, 6, 2008.

\bibitem{MR2559677}
Christiane Quesne.
\newblock Solvable rational potentials and exceptional orthogonal polynomials in supersymmetric quantum mechanics.
\newblock {\em SIGMA Symmetry Integrability Geom. Methods Appl.}, 5:Paper 084, 24, 2009.

\bibitem{Slater}
L. J. Slater, \emph{Generalized hypergeometric functions}, Cambridge University Press, Cambridge, 1966.


\bibitem{Szego}
G. Szeg\"{o}, \emph{Orthogonal Polynomials}. Providence (RI): American Mathematical Society, 1975.


\end{thebibliography}
\end{document}